\documentclass{amsart}

\newtheorem{theorem}{Theorem}[section]
\newtheorem{lemma}[theorem]{Lemma}
\newtheorem{proposition}[theorem]{Proposition}
\newtheorem{corollary}[theorem]{Corollary}

\theoremstyle{definition}

\newtheorem{conjecture}[theorem]{Conjecture}

\theoremstyle{remark}
\newtheorem{remark}[theorem]{Remark}

\newcommand{\R}{\mathbb{R}}
\newcommand{\C}{\mathbb{C}}
\newcommand{\ad}{\mathrm{ad}}
\newcommand{\Ad}{\mathrm{Ad}}
\newcommand{\rank}{\mathrm{rank}}
\newcommand{\SO}{\mathrm{SO}}
\newcommand{\GL}{\mathrm{GL}}

\newcommand{\SL}{\mathrm{SL}}

\newcommand{\Lie}{\mathrm{Lie}}
\newcommand{\Aut}{\mathrm{Aut}}
\newcommand{\Kill}{\mathrm{Kill}}
\newcommand{\Iso}{\mathrm{Iso}}

\newcommand{\OO}{\mathcal{O}}

\newcommand{\GG}{\mathcal{G}}

\newcommand{\g}{\mathfrak{g}}
\newcommand{\h}{\mathfrak{h}}
\newcommand{\gl}{\mathfrak{gl}}
\newcommand{\so}{\mathfrak{so}}
\newcommand{\slinear}{\mathfrak{sl}}
\newcommand{\s}{\mathfrak{s}}

\title[Actions and transverse structures: Integrable normal]{Isometric actions
of simple groups and transverse structures: The integrable normal case}
\author{Raul Quiroga-Barranco}
\address{Centro de Investigaci\'on en Matem\'aticas,
Apartado Postal 402, 36000, Guanajuato, Guanajuato, M\'exico.}
\email{quiroga@cimat.mx}

\dedicatory{With gratitude to Robert J.~Zimmer on the occasion of his 60th
birthday.}

\thanks{This work was supported by SNI-Mexico and the grants Conacyt 44620 and
Concyteg 07-02-K662-091.}

\begin{document}

\maketitle

\begin{abstract}
For actions with a dense orbit of a connected noncompact simple Lie group $G$,
we obtain some global rigidity results when the actions preserve certain
geometric structures. In particular, we prove that for a $G$-action to be
equivalent to one on a space of the form $(G\times K\backslash H)/\Gamma$, it is
necessary and sufficient for the $G$-action to preserve a pseudo-Riemannian
metric and a transverse Riemannian metric to the orbits. A similar result proves
that the $G$-actions on spaces of the form $(G\times H)/\Gamma$ are
characterized by preserving transverse parallelisms. By relating our techniques
to the notion of the algebraic hull of an action, we obtain infinitesimal Lie
algebra structures on certain geometric manifolds acted upon by $G$.
\end{abstract}

\section{Introduction}
\label{section-intro}
In the rest of this work, we let $G$ be a connected noncompact simple Lie
group with Lie algebra $\g$ and $M$ a smooth connected manifold acted
upon smoothly by $G$. There are several examples of such actions that preserve a
finite volume.
Some of the most interesting are obtained from Lie group homomorphisms $G
\hookrightarrow L$ where $L$ is a connected Lie group that admits a lattice
$\Gamma$. The relevant $G$-action is then given on the double coset
$K\backslash L/\Gamma$, where $K$ is some compact subgroup that centralizes $G$.
Robert Zimmer proposed in \cite{Zimmer-prog} to study the finite measure
preserving $G$-actions on $M$ and determine to what extent these are given by
such double cosets.

To understand how to tackle Zimmer's program, it has been proved very useful to
consider $G$-actions preserving a suitable geometric structure (see
\cite{Gromov} and \cite{Zimmer-geometric}). This is a natural condition since
the above double coset examples carry a $G$-invariant pseudo-Riemannian metric
when $L$ is semisimple. Such metric comes from the bi-invariant metric on $L$
obtained from the Killing form of its Lie algebra.

The properties of the $G$-orbits of finite measure preserving $G$-actions are
reasonably well known. For example, such $G$-actions are known to be everywhere
locally free when they preserve suitable pseudo-Riemannian metrics (see
\cite{Gromov}, \cite{Szaro} and \cite{Zeghib}). Also in such case, the metric
restricted to the orbits can be
described precisely (see Lemma~\ref{lemma-metric-on-orbits}). However, there is
still a lack of knowledge of the properties of a manifold, acted upon by
$G$, in the transverse direction to the orbits. 

In this work, we propose to study the $G$-actions on $M$ by emphasizing the
need to understand the properties of the transverse to the $G$-orbits. For
this, we will be dealing with $G$-actions that have a dense orbit and preserve
a finite volume pseudo-Riemannian metric. As observed above, in this case the
results from \cite{Szaro} show that the $G$-orbits define a foliation, which
from now on we will denote with $\OO$. By assuming that the orbits are
nondegenerate for the pseudo-Riemannian metric on $M$, we can consider the
normal bundle $T\OO^\perp$ to the orbits as realizing the transverse direction
in $M$. With this respect, we obtain the following structure theorem for
$G$-actions on $M$ when this normal bundle is integrable.

\begin{theorem}
\label{thm-int-normal}
Suppose that the $G$-action on $M$ has a dense orbit and preserves a finite
volume complete pseudo-Riemannian metric. If the $G$-orbits in $M$ are
nondegenerate and the normal bundle to the orbits $T\OO^\perp$ is integrable,
then there exist:
\begin{enumerate}
\item an isometric finite covering map $\widehat{M} \rightarrow M$ to which the
$G$-action lifts,

\item a simply connected complete pseudo-Riemannian manifold $\widetilde{N}$,
and

\item a discrete subgroup $\Gamma \subset G \times \Iso(\widetilde{N})$,
\end{enumerate}
such that $\widehat{M}$ is $G$-equivariantly isometric to $(G\times
\widetilde{N})/\Gamma$.
\end{theorem}

We observe that in the double coset examples $K\backslash L/\Gamma$ as above
with $L$ semisimple and with the metric coming from the Killing form, the
$G$-orbits are always nondegenerate. Besides that, we prove that for a general
$G$-action on $M$, the orbits are always nondegenerate when $\dim(M) < 2\dim(G)$
(see Lemma~\ref{lemma-nondeg}). Also, by developing criteria for the normal
bundle $T\OO^\perp$ to be integrable, we obtain some results where the
conclusion of Theorem~\ref{thm-int-normal} holds. For example,
Corollary~\ref{cor-int-normal-rank} ensures such conclusion when $G$ has high
enough real rank. Also, Corollary~\ref{cor-int-normal-reps} does the same for
manifolds $M$ whose dimensions have a suitable bound in terms of $G$. In this
last result we can even dispense with the assumption of having nondegenerate
orbits by applying Lemma~\ref{lemma-nondeg}.

Without assuming in Theorem 1.1 that the $G$-action has a dense orbit it is
still possible to draw some conclusions. More specifically, if we assume the
rest of the hypotheses, a description of the universal covering space of the
manifold as a warped product is obtained. This sort of result already appears in
\cite{Hernandez}.

For a $G$-action on $M$ preserving a pseudo-Riemannian metric, in
\cite{Quiroga-Annals} we considered a certain geometric condition between the
metrics on $G$ and $M$ (the former for a bi-invariant metric) from which we
concluded that $M$ is, up to a finite covering space, a double coset of the form
$(G\times K\backslash H)/\Gamma$. In other words, a double coset as above
where $G$ appears as a factor in $L$. One of the steps used in
\cite{Quiroga-Annals} to achieve this was to prove that the normal bundle
$T\OO^\perp$ is Riemannian. Considering the relevance we are giving to the
transverse to the orbits, it is natural to determine the properties of the
$G$-actions that preserve a transverse Riemannian structure. In this context, we
obtain the following result which proves that, up to a finite covering, the
double cosets of the form $(G\times K\backslash H)/\Gamma$ are characterized as
those isometric $G$-actions that preserve a transverse Riemannian structure on
the foliation $\OO$.
We recall that a semisimple Lie group is isotypic if the complexification of its
Lie algebra is isomorphic to a sum of identical simple complex ideals.

\begin{theorem}\label{theorem-trans-Riem-equiv}
If $G$ is a connected noncompact simple Lie group acting faithfully on a
compact manifold $M$, then the following conditions are equivalent:
\begin{enumerate}
\item There is a $G$-equivariant finite covering map $(G\times K\backslash
H)/\Gamma \rightarrow M$ where $H$ is a connected Lie group with a compact
subgroup $K$ and $\Gamma \subset G \times H$ is a discrete cocompact subgroup
such that $G\Gamma$ is dense in $G \times H$.

\item There is a finite covering map $\widehat{M} \rightarrow M$ for which the
$G$-action on $M$ lifts to a $G$-action on $\widehat{M}$ with a dense orbit and
preserving:
\begin{itemize}
\item a pseudo-Riemannian metric for which the orbits are nondegenerate,
\item a transverse Riemannian structure for the foliation $\OO$ by $G$-orbits.
\end{itemize}
\end{enumerate}
Furthermore, if $G$ has finite center and real rank at least $2$, then we can
assume in (1) of the above equivalence that $G \times H$ is semisimple isotypic
with finite center and that $\Gamma$ is an irreducible lattice.
\end{theorem}

Based on the previous result, we prove in Theorem~\ref{theorem-Lorentzian} that
the $G$-actions on $M$ preserving a Lorentzian metric are, up to a finite
covering, those given by double cosets $(G \times K\backslash H)/\Gamma$
with $G$ locally isomorphic to $\SL(2,\R)$. We observe that this result
improves a similar one found in \cite{Gromov} (see also the Introduction of
\cite{Nadine-Annals}).

Continuing with our study of transverse structures, we next consider isometric
$G$-actions preserving a transverse parallelism for the foliation $\OO$. We
prove that, up to a finite covering, such actions characterize the double cosets
of the form $(G\times H)/\Gamma$.

\begin{theorem}\label{theorem-trans-parallelism-equiv}
If $G$ is a connected noncompact simple Lie group acting faithfully on a
compact manifold $M$, then the following conditions are equivalent:
\begin{enumerate}
\item There is a $G$-equivariant finite covering map $(G\times H)/\Gamma
\rightarrow M$ where $H$ is a connected Lie group and $\Gamma \subset G \times
H$ is a discrete cocompact subgroup such that $G\Gamma$ is dense in $G \times
H$.

\item There is a finite covering map $\widehat{M} \rightarrow M$ for which the
$G$-action on $M$ lifts to a $G$-action on $\widehat{M}$ with a dense orbit and
preserving:
\begin{itemize}
\item a pseudo-Riemannian metric for which the orbits are nondegenerate,
\item a transverse parallelism for the foliation $\OO$ by $G$-orbits.
\end{itemize}
\end{enumerate}
Furthermore, if $G$ has finite center and real rank at least $2$, then we can
assume in (1) of the above equivalence that $G \times H$ is semisimple isotypic
with finite center and that $\Gamma$ is an irreducible lattice.
\end{theorem}

The notion of the algebraic hull of an action on a bundle, introduced by
Zimmer, is a fundamental tool to understand $G$-actions as they relate to
geometric structures. We recall that the algebraic hull is the smallest
algebraic subgroup for which there is a measurable $G$-invariant reduction of
the bundle being acted upon (see \cite{Zimmer-alghull}). For our
setup, where we are interested in the transverse to the orbits, it is then
natural to consider the algebraic hull of the $G$-action on the bundle
$L(T\OO^\perp)$ for an isometric $G$-action with nondegenerate orbits. Since
a $G$-action of this sort preserves a pseudo-Riemannian metric on $T\OO^\perp$,
the algebraic hull for $L(T\OO^\perp)$ is in this case a subgroup of
$\mathrm{O}(p,q)$, for some $(p,q)$. The next result shows that for weakly
irreducible manifolds and when the algebraic hull for $L(T\OO^\perp)$ is the
largest possible for this setup, the manifold being acted upon has
infinitesimally at some point the structure of a specific Lie algebra that
contains $\g$. By the last claim in the statement, such Lie algebra structure is
nontrivially linked to the geometry of the manifold. For the explicit
description of the Lie algebra structure obtained in this result we refer to
Lemma~\ref{lemma-Omega-curvature} and Theorem~\ref{theorem-alghull-TOperp}.

We recall that a pseudo-Riemannian manifold is weakly irreducible if the tangent
space at some (and hence any) point has no proper nondegenerate subspaces
invariant under the restricted holonomy group at that point. Also recall that
every simple Lie group with a bi-invariant metric is weakly irreducible.

\begin{theorem}\label{theorem-alghull-TOperp-irred}
Suppose that $G$ has finite center and real rank at least $2$, and that the
$G$-action on $M$ preserves a finite volume complete pseudo-Riemannian metric.
Also assume that $G$ acts ergodically on $M$ and that the foliation by
$G$-orbits is nondegenerate. Denote with $L$ the algebraic hull for the
$G$-action on the bundle $L(T\OO^\perp)$ and with $\mathfrak{l}$ its Lie
algebra. In particular, there is an embedding of Lie algebras $\mathfrak{l}
\hookrightarrow \so(p,q)$, where $(p,q)$ is the signature of the metric of $M$
restricted to $T\OO^\perp$. If this embedding is surjective and $M$ is weakly
irreducible, then the following holds:
\begin{enumerate}
\item $G$ is locally isomorphic to $\SO_0(p,q)$ and $\dim(M) = (p+q)(p+q+1)/2$,

\item for some $x \in M$ the tangent space $T_x M$
admits a Lie algebra structure isomorphic to either $\so(p,q+1)$ or
$\so(p+1,q)$ such that $T_x\OO$ is a Lie subalgebra isomorphic to $\g$.
\end{enumerate}
Furthermore, with respect to the representation of Lemma~\ref{lemma-lambda},
there is a Lie algebra of local Killing fields vanishing at $x$ which is
isomorphic to $\g$ and acts nontrivially on $T_x M$ by derivations
of the Lie algebra structure given in (2).
\end{theorem}

With the results developed in this work, we try to show the importance of
considering transverse geometric structures to understand the actions of
noncompact simple Lie groups. In fact, we believe that some form of the
following conjecture could provide a geometric characterization of the double
coset examples of $G$-actions mentioned above.

\begin{conjecture}
Consider the double coset $G$-spaces of the form $K\backslash L/\Gamma$, where
$L$ is a semisimple Lie group with an irreducible lattice $\Gamma$ and a compact
subgroup $K$; the $G$-action is then induced by a nontrivial homomorphism $G
\rightarrow L$ whose image centralizes $K$. Then, a $G$-action on a manifold $M$
is equivalent to such a double coset $G$-action for some $L,\Gamma,K$ if and
only if the $G$-action on $M$:
\begin{itemize}
\item has a dense orbit,
\item preserves a pseudo-Riemannian metric, and
\item preserves a transverse geometric structure to the orbits suitably related
to the geometry of $GK\backslash L$.
\end{itemize}
\end{conjecture}

Note that for $L = G\times H$ and $K$ a compact subgroup of $H$, we obtain
quotients $GK\backslash L = K\backslash H$ and $G\backslash L = H$ which
naturally carry a Riemannian metric and a parallelism, respectively. Thus
Theorems~\ref{theorem-trans-Riem-equiv} and
\ref{theorem-trans-parallelism-equiv} verify the conjecture in these cases.

One of the main tools used to obtain our results is Gromov's machinery on
geometric $G$-actions (see \cite{GCT}, \cite{Gromov}, \cite{Quiroga-Annals} and
\cite{Zimmer-geometric}). Such machinery ensures the existence of large families
of local Killing fields. We develop these techniques in
Section~\ref{section-Killing-fields} making emphasis on the fact that the
Killing fields thus obtained yield $\g$-module structures on the tangent spaces
to $M$; with this respect, our main result is Proposition~\ref{prop-g(x)} which
is due to Gromov~\cite{Gromov} in the analytic/compact case (see also
\cite{Zimmer-geometric}) and was extended to the smooth/finite volume case in
\cite{GCT} and \cite{Quiroga-Annals}.
In Section~\ref{section-int-normal}, we prove that the latter impose
restrictions strong enough to guarantee the integrability of the normal bundle
under suitable conditions (Corollaries~\ref{cor-int-normal-rank} to
\ref{cor-int-compact-Riemannian}). We observe that Theorem~\ref{thm-int-normal}
and its consequences obtained in Section~\ref{section-int-normal} are in fact
extensions of results obtained in \cite{Gromov}: Theorem~5.3.E in page~129 of
\cite{Gromov} states, under the assumptions of our
Corollary~\ref{cor-int-normal-rank}, that $M$ has a covering $G$-equivariantly
diffeomorphic to $G \times N$ for some $N$. Note, however, that the result in
\cite{Gromov} yields only a topological covering and does not further describes
the covering group.

Theorem~\ref{theorem-trans-Riem-equiv} is obtained in
Section~\ref{section-Riem} by applying the main results and arguments from
\cite{Quiroga-Annals}. As mentioned above, this yields in
Section~\ref{section-Lorentz} our characterization of Lorentzian manifolds
acted upon with a dense orbit by a simple noncompact Lie group. In
Section~\ref{section-parallel} we establish the characterization of $G$-spaces
of the form $(G\times H)/\Gamma$ provided by
Theorem~\ref{theorem-trans-parallelism-equiv}; for this, one of the main
ingredients is given by Theorem~\ref{thm-int-normal}. We observe that the
arguments used in Section~\ref{section-parallel} are based on those found in
\cite{Quiroga-AJM}.

To obtain Theorem~\ref{theorem-alghull-TOperp-irred} in
Section~\ref{section-alghull}, an application of Theorem~\ref{thm-int-normal} is
used again. But now, the notion of the algebraic hull and the deep result about
it found in \cite{Zimmer-alghull} are also fundamental. Ultimately, this is
somewhat to be expected, since a computation of the algebraic hull for frame
bundles over $M$ is in fact an important step to build Gromov's machinery for
geometric $G$-actions on $M$. This is evident in our proof of
Proposition~\ref{prop-g(x)}.

We want to observe that the conclusion of
Theorem~\ref{theorem-alghull-TOperp-irred} can be paraphrased by saying that
the manifold $M$ has infinitesimally the structure of either $\SO(p,q+1)$ or
$\SO(p+1,q)$. It remains the question as to whether or not $M$ can be related
more precisely to either of these groups. This problem will be pursued
elsewhere (see \cite{Gestur-Quiroga}) by requiring some additional conditions.

The author wishes to thank Uri Bader and Shmuel Weinberger for some fruitful
conversations.

\section{Killing fields on manifolds with a simple group of isometries}
\label{section-Killing-fields}
By Theorem 4.17 from \cite{Szaro}, if the $G$-action on $M$ has a dense orbit
and preserves a finite volume pseudo-Riemannian metric, then the action is
locally free and so
the orbits define a foliation that we have agreed to denote with $\OO$. In this
case, it is well known that the bundle $T\OO$ tangent to the foliation $\OO$ is
a trivial vector bundle isomorphic to $M \times \g$, under the isomorphism given
by:
\begin{align*}
    M \times \g &\rightarrow T\OO \\
        (x, X) &\mapsto X^*_x.
\end{align*}
For every $x \in M$, this induces an isomorphism between
the fiber $T_x\OO$ and $\g$, which we will refer to as the natural isomorphism.
Furthermore, if we consider on $M\times \g$ the product $G$-action, where the
$G$-action on $\g$ is the adjoint one, then such isomorphism is
$G$-equivariant.  More precisely, we have:
$$
	dg(X^*) = \Ad(g)(X)^*
$$
for every $g \in G$ and $X \in \g$. 
Note that for $X$ in the Lie algebra of a group acting on a manifold, we denote
with $X^*$ the vector field on the manifold whose one-parameter group of
diffeomorphisms is given by $(\exp(tX))_t$ through the action on the manifold.

For any given pseudo-Riemannian manifold $N$, we will denote with $\Kill(N,x)$
the Lie algebra of germs at $x$ of local Killing vector fields defined in a
neighborhood of $x$, and with $\Kill_0(N,x)$ we will denote the Lie subalgebra
consisting of those germs that vanish at $x$. The following result is a
consequence of the Jacobi identity and the fact that the Lie derivative of
a metric with respect to its Killing fields vanishes. In the rest of this work,
for a vector space $W$ with a nondegenerate symmetric bilinear form, we will
denote with $\so(W)$ the Lie algebra of linear maps on $W$ that are
skew-symmetric with respect to the bilinear form.

\begin{lemma}\label{lemma-lambda}
Let $N$ be a pseudo-Riemannian manifold and $x\in N$. Then, the map:
\begin{align*}
    \lambda_x : \Kill_0(N,x) &\rightarrow \so(T_x N) \\
        \lambda_x(Z)(v) &= [Z,V]_x,
\end{align*}
where $V$ is any vector field such that $V_x = v$, is a well defined
homomorphism of Lie algebras.
\end{lemma}

From now on, for a given point $x$ of a pseudo-Riemannian manifold, the map
$\lambda_x$ will denote the homomorphism from the previous lemma.

For the proof of our next result we will present some facts about infinitesimal
automorphisms and Killing fields, and we refer to \cite{GCT} for further
details. The tangent bundle of order $k$ of a manifold $N$ is the smooth bundle
$T^{(k)}N$ whose fiber at $x$ is the space $T_x^{(k)}N$ of $(k-1)$-jets at $x$
of vector fields of $N$. For every (local) diffeomorphism $\varphi : N_1
\rightarrow N_2$ mapping $x_1$ to $x_2$ we have a linear isomorphism:
\begin{align*}
	T_{x_1}^{(k)}N_1 &\rightarrow T_{x_2}^{(k)}N_2  \\
		j^{k-1}_{x_1}(X) &\mapsto j^{k-1}_{x_2}(d\varphi(X))
\end{align*}
that depends only on the jet $j^k_{x_1}(\varphi)$. 
For $N_1=N_2=N$, $x_1=x_2=x$ this yields the group
$D^{(k)}_x(N)$ of $k$-jets at $x$ of local diffeomorphisms fixing $x$,
whose group structure is induced by the composition of maps. In the case when $N
= \R^n$ and $x = 0$, we will denote this group with $\GL^{(k)}(n)$. 
We also recall that the $k$-th order frame bundle over $N$ is the smooth bundle
$L^{(k)}(N)$ that consists of the $k$-jets at $0$ of local diffeomorphisms $\R^n
\rightarrow N$; any such a jet $j^k_0(\varphi)$ thus defines a linear
isomorphism $T_0^{(k)}\R^n \rightarrow T_{\varphi(0)}^{(k)}N$. The structure
group of $L^{(k)}(N)$ is $\GL^{(k)}(n)$.
With respect to vector fields, we denote with $\mathcal{D}^{(k)}_x(N)$ the space
of $k$-jets at $x$ of vector fields vanishing at $x$, and we use the special
notation $\gl^{(k)}(n)$ when $N = \R^n$ and $x = 0$. 
If $N$ carries a pseudo-Riemannian metric, for every $x \in N$, we will denote
with $\Aut^k(N,x)$ the subgroup of $D^{(k)}_x(N)$ consisting of those $k$-jets
that preserve the metric up to order $k$ at $x$. 
Correspondingly, for vector fields we denote with $\Kill^k(N,x)$ the space of
$k$-jets at $x$ of vector fields on $N$ that preserve the metric up to order $k$
at $x$; the subspace of those $k$-jets whose
$0$-jet vanishes is denoted with $\Kill_0^k(N,x)$. The next result provides a
natural representation of $D^{(k)}_x(N)$ from which the Lie algebras of this
group and of $\Aut^k(N,x)$ are described in terms of $\mathcal{D}^{(k)}_x(N)$
and $\Kill_0^k(N,x)$, respectively; the proof is elementary, but it is detailed
in Section~2 and 4 of \cite{GCT}.

\begin{lemma}\label{lemma-jets}
For a smooth manifold $N$ and any given point $x \in N$ the following
properties hold for every $k\geq 1$:
\begin{enumerate}
 \item The map 
\begin{align*}
	\Theta_x : D^{(k)}_x(N) &\rightarrow \GL(T_x^{(k)}N)  \\
 	\Theta_x(j^k_x(\varphi))(j^{k-1}_x(X)) &= j^{k-1}_x(d\varphi(X))
\end{align*}
is a Lie group monomorphism. 

\item The assignment $[j^k_x(X),j^k_x(Y)]^k = -j^k_x([X,Y])$
yields a well defined Lie algebra structure on $\mathcal{D}^{(k)}_x(N)$.

\item The map:
\begin{align*}
	\theta_x : \mathcal{D}^{(k)}_x(N) &\rightarrow \gl(T_x^{(k)}N)  \\
 	\theta_x(j^k_x(X))(j^{k-1}_x(Y)) &= -j^{k-1}_x([X,Y])
\end{align*}
is a Lie algebra monomorphism for the Lie algebra structure on
$\mathcal{D}^{(k)}_x(N)$ given by $[\cdot,\cdot]^k$. Furthermore,
$\theta_x(\mathcal{D}^{(k)}_x(N)) = \Lie(\Theta_x(D^{(k)}_x(N)))$.

\item If $N$ has a pseudo-Riemannian metric, then we have
$\theta_x(\Kill_0^k(N,x)) = \Lie(\Theta_x(\Aut^k(N,x)))$.
\end{enumerate}
In particular, with respect to the homomorphisms $\Theta_x$ and $\theta_x$, the
Lie algebra of $\Aut^k(N,x)$ is realized by $\Kill_0^k(N,x)$ with the Lie
algebra structure given by $[\cdot,\cdot]^k$.
\end{lemma}

The following result is due to Gromov in the analytic case (see
\cite{Gromov}) and it was extended to the smooth case in~\cite{Quiroga-Annals}.
We present here a fairly detailed proof based on the results from \cite{GCT}.
In congruence with our notation for $M$, in the rest of this 
work we will use $\OO$ to denote the foliation by $\widetilde{G}$-orbits in
$\widehat{M}$ for the lifted $\widetilde{G}$-action on every covering space
$\widehat{M}$ of $M$; this will be so for any covering whether finite or not.

\begin{proposition}\label{prop-g(x)}
Suppose that the $G$-action on $M$ has a dense orbit and preserves a finite
volume pseudo-Riemannian metric. Then, there exists a dense subset $S\subset
\widetilde{M}$ such that for every $x \in S$ the following properties are
satisfied.
\begin{enumerate}
\item There is a homomorphism of Lie algebras $\rho_x : \g \rightarrow
\Kill(\widetilde{M},x)$ which is an isomorphism onto its image
$\rho_x(\g) = \g(x)$.

\item $\g(x)\subset \Kill_0(\widetilde{M},x)$, i.e.~every element of $\g(x)$
vanishes at $x$.

\item For every $X,Y \in \g$ we have:
$$
    [\rho_x(X),Y^*] = [X,Y]^* = -[X^*,Y^*],
$$
in a neighborhood of $x$. In particular, the elements in $\g(x)$ and their
corresponding local flows preserve both $\OO$ and $T\OO^\perp$ in a
neighborhood of $x$.

\item The homomorphism of Lie algebras $\lambda_x\circ\rho_x : \g \rightarrow
\so(T_x \widetilde{M})$ induces a $\g$-module structure on $T_x \widetilde{M}$
for which the subspaces $T_x \OO$ and $T_x \OO^\perp$ are $\g$-submodules.
\end{enumerate}
\end{proposition}
\begin{proof}
Since the proof builds on the notions and results found in \cite{GCT}, we
will mostly follow its notation. We will be careful to define the objects
considered but we refer to \cite{GCT} for further details.

For every $k \geq 1$, let us denote with $\sigma^k : L^{(k)}(M) \rightarrow Q_k$
the $\GL^{(k)}(n)$-equivariant map that defines the $k$-th order extension of
the geometric structure defined by the pseudo-Riemannian metric on $M$.

Consider the set:
$$
	\GG^k = \{ j^{k-1}_x(X^*) : X \in \g, x\in M\},
$$
which, by the local freeness of the $G$-action, is a smooth subbundle of
$T^{(k)}M$. In fact, we have $\GG = T\OO$, and as with this bundle there is a
trivialization given by:
\begin{align*}
    M \times \g &\rightarrow \GG^k \\
        (x, X) &\mapsto j^{k-1}_x(X^*).
\end{align*}
The corresponding trivialization of the frame bundle of $\GG^k$ is given by:
\begin{align*}
    M \times \GL(\g) &\rightarrow L(\GG^k) \\
        (x, A) &\mapsto A_x.
\end{align*}
where $A_x(X) = j^{k-1}_x((AX)^*)$. Note that we have taken $\g$ as the
standard fiber of the bundle $\GG^k$.

Choose a subspace $\GG_0$ of $T_0^{(k)}\R^n$ isomorphic to $\g$. We will fix
such subspace as well as an isomorphism with $\g$ through which we will
identify these two spaces. Let us now consider:
$$
	L^{(k)}(M,\GG^k) = \{ \alpha \in L^{(k)}(M) : 
		\alpha(\GG_0) = \GG^k_x \mbox{ if } \alpha \in L^{(k)}(M)_x \}
$$
which is a smooth reduction of $L^{(k)}(M)$ to the subgroup of $\GL^{(k)}(n)$
that preserves the subspace $\GG_0$; we will denote such subgroup with
$\GL^{(k)}(n,\GG_0)$. Recall from the remarks preceding Lemma~\ref{lemma-jets}
that for every $j^k_0(\varphi) \in L^{(k)}(M)$ we obtain a linear isomorphism:
\begin{align*}
	T_0^{(k)} \R^n & \rightarrow T_{\varphi(0)}^{(k)} M  \\
		j^{k-1}_0(X) &\mapsto j^{k-1}_{\varphi(0)}(d\varphi(X)).
\end{align*}
In particular, if we let:
\begin{align*}
	f_k : L^{(k)}(M,\GG^k) &\rightarrow L(\GG^k)  \\
		j^k_0(\varphi) &\mapsto j^k_0(\varphi)|_{\GG_0}.
\end{align*}
then, by the identification between $\GG_0$ and $\g$, we can consider $f_k$ as a
well-defined smooth principal bundle morphism that covers the
identity. The associated homomorphism of structure groups for $f_k$ is given by:
\begin{align*}
	\pi_k : \GL^{(k)}(n,\GG_0) &\rightarrow \GL(\g)  \\
		j^k_0(\varphi) &\mapsto j^k_0(\varphi)|_{\GG_0},
\end{align*}
which is clearly surjective. Note that we have used again our identification
between $\g$ and $\GG_0$.

Fix $\mu$ an arbitrary ergodic component for the $G$-action on $M$ for the
pseudo-Riemannian volume. Then, there is a measurable reduction $P$ of
$L^{(k)}(M,\GG^k)$ so that $\sigma^k(P)$ is ($\mu$-a.e.~over $M$) a single point
$q_0 \in Q_k$. Furthermore, the structure group of $P$ is the subgroup of
$\GL^{(k)}(n,\GG_0)$ that stabilizes $q_0$, and in particular it is
algebraic. This claim is a consequence of the fact that the
$\GL^{(k)}(n,\GG_0)$-action on $Q_k$ is algebraic, which in turn follows from
the fact that a pseudo-Riemannian metric is a geometric structure of algebraic
type; we refer to Section~4 and the proof of Proposition~8.4 of \cite{GCT} for
further details.

On the other hand, since $\pi_k$ is a surjection and $f_k$ is $G$-equivariant,
by Proposition~8.2 of \cite{GCT}, there exist reductions $Q_1$ and $Q_2$ of
$L^{(k)}(M,\GG^k)$ and $L(\GG^k)$, respectively, to subgroups $L_1 \subset
\GL^{(k)}(n,\GG_0)$ and $\overline{\Ad(G)}^Z \subset \GL(\g)$, such that
$f_k(Q_1) \subset Q_2$ ($\mu$-a.e.~over $M$) and such that $\pi_k(L_1)$ is a
finite index subgroup of $\overline{\Ad(G)}^Z$. 
Here $L_1$ can be chosen to be the algebraic hull of $L^{(k)}(M,\GG^k)$ for the
$G$-action on $M$ with respect to the ergodic measure $\mu$. 
This claim uses the well known fact that $\overline{\Ad(G)}^Z$ is the
algebraic hull of $M \times \GL(\g)$ for the product action. It is not
difficult to see that this can be chosen so that $Q_2 = M
\times \overline{\Ad(G)}^Z$ ($\mu$-a.e.~over $M$) with respect to
the above identification $M \times \GL(\g) \cong L(\GG^k)$.
We can also assume that $Q_1 \subset P$, $\mu$-a.e.~over $M$, because
the reduction $Q_1$ is the smallest one to an algebraic subgroup.

The above discussion ensures that for $\mu$-a.e.~$x \in M$, we have the
following relations:
\begin{align*}
	f_k((Q_1)_x) &\subset (Q_2)_x = \{x\}\times\overline{\Ad(G)}^Z  \\
	(Q_1)_x &\subset (P)_x \subset L^{(k)}(M,\GG^k)_x  \\
	\sigma^k((P)_x) &= \{q_0\} 
\end{align*}
Let us now fix a point $x$ such that these conditions hold. 
Choose $\alpha_x \in (Q_1)_x$ and let $f_k(\alpha_x) = (x,k_x)$ where $k_x \in
\overline{\Ad(G)}^Z$. Since $\pi_k$ is surjective, there exists $\widehat{k}_x
\in \GL^{(k)}(n,\GG_0)$ such that $\pi_k(\widehat{k}_x) = k_x$. In particular,
by the $\pi_k$-equivariance of $f_k$ we have $f_k(\alpha_x \widehat{k}_x^{-1}) =
(x,e)$. We also have by the same reason:
$$
	f_k(\alpha_x g \widehat{k}_x^{-1}) 
		= f_k(\alpha_x \widehat{k}_x^{-1} \widehat{k}_x g 
			\widehat{k}_x^{-1})
		= (x,k_x \pi_k(g) k_x^{-1}),
$$
for every $g \in L_1$. Also, the inclusion $(Q_1)_x \subset
L^{(k)}(M,\GG^k)_x$ implies that, for every $g \in L_1$, the $k$-jets of
diffeomorphisms $\alpha_x\widehat{k}_x^{-1}, \alpha_x g\widehat{k}_x^{-1}$
considered as linear isomorphisms $T_0^{(k)}\R^n \rightarrow T_x^{(k)}M$ map
$\GG_0$ onto $\GG^k_x$. Hence, from the definition of $f_k$ it follows that
$\alpha_x g \alpha_x^{-1} = (\alpha_x g\widehat{k}_x^{-1})
(\alpha_x\widehat{k}_x^{-1})^{-1}$ is a $k$-jet of local diffeomorphism 
of $M$ at $x$ whose associated isomorphism $T_x^{(k)}M \rightarrow T_x^{(k)}M$
maps $\GG^k_x$ onto itself by the assignment:
$$
	j^{k-1}_x(X^*) \mapsto j^{k-1}_x((k_x \pi_k(g) k_x^{-1} X)^*).
$$
for which we have used the above trivialization $M\times \GL(\g) \cong
L(\GG^k)$. Since $\pi_k(L_1)$ has finite index in $\overline{\Ad(G)}^Z$ it
contains the identity component $\Ad(G)$, and because $k_x \in
\overline{\Ad(G)}^Z$ the group $k_x \pi_k(L_1) k_x^{-1}$ also contains $\Ad(G)$.
It follows that $\alpha_xL_1\alpha_x^{-1}$ is a subgroup of
$D^{(k)}_x(M)$ for which the homomorphism from Lemma~\ref{lemma-jets}(1)
induces a homomorphism:
\begin{align*}
	H_x : \alpha_xL_1\alpha_x^{-1} &\rightarrow \GL(\GG^k_x)  \\
	\alpha_x g \alpha_x^{-1} &\mapsto 
		\Theta_x(\alpha_x g \alpha_x^{-1})|_{\GG^k_x}
\end{align*}
whose image contains $\Ad(G) \subset \GL(\g)$ with respect to the
identification between $\g$ and $\GG^k_x$ given by the above isomorphism
$M\times \g \cong \GG^k$. This implies that the induced Lie algebra
homomorphism:
$$
	h_x : \Lie(\alpha_xL_1\alpha_x^{-1}) 
		\rightarrow \gl(\GG^k_x)
$$
has image $\ad(\g)$, again with respect to the referred identification between
$\g$ and $\GG^k_x$.

On the other hand, we have for every $g \in L_1$:
$$
	\sigma^k((\alpha_x g \alpha_x^{-1})\alpha_x) =
		\sigma^k(\alpha_x g) = \sigma^k(\alpha_x)
$$
because $\sigma^k((Q_1)_x) \subset \sigma^k((P)_x) = \{q_0\}$ is a single point.
But this identity proves that every such $k$-jet $\alpha_x g \alpha_x^{-1}$
preserves the pseudo-Riemannian metric up order $k$ (see \cite{GCT}). In other
words, $\alpha_xL_1\alpha_x^{-1}$ is a subgroup of $\Aut^k(M,x)$ and by
Lemma~\ref{lemma-jets} we also have that $\Lie(\alpha_xL_1\alpha_x^{-1})$ is a
Lie subalgebra of $\Kill_0^k(M,x)$.

From the above remarks, it follows that there is a Lie algebra homomorphism
$\widehat{\rho}^k_x : \g \rightarrow \Kill^k_0(M,x)$ such that:
\begin{equation*}
	\theta_x(\widehat{\rho}^k_x(X))(j^{k-1}_x(Y^*)) 
		= j^{k-1}_x([X,Y]^*) \quad \mbox{ for every } 
		X,Y \in \g. \tag{*}
\end{equation*}

For $k$ fixed, the existence of the homomorphism $\widehat{\rho}^k_x$ has been
established for $\mu$-a.e.~$x \in M$, where $\mu$ is an arbitrary ergodic
component of the pseudo-Riemannian volume of $M$. Thus, for $k$ fixed, it
follows that the homomorphism $\widehat{\rho}^k_x$ exists for every $x \in S_k$,
where $S_k$ is some subset of $M$ which is conull with respect to the
pseudo-Riemannian volume of $M$. 
Finally, if we let $S_0 = \cap_{k=1}^\infty S_k$, then $S_0$ is conull
with respect to the pseudo-Riemannian volume and for every $x \in S_0$ and every
$k \geq 1$ there exist a homomorphism $\widehat{\rho}^k_x : \g \rightarrow
\Kill^k_0(M,x)$ satisfying~(*).

In~\cite{Nomizu} the notion of $\mathfrak{k}$-regular point for a metric in a
manifold is
introduced. Such regular points satisfy two key properties relevant to our
discussion. First, the set of regular points $U$ of $M$ is an open dense subset.
Second, for $x\in U$ there is some integer $k(x) \geq 1$ so that, for $k \geq
k(x)$, every element of $\Kill^k_0(M,x)$ extends uniquely to an element of
$\Kill_0(M,x)$.  The first property is found in \cite{Nomizu} and the second one
is proved in \cite{GCT}, both just using smoothness. Note that the results in
\cite{Nomizu} are stated for Riemannian manifolds but, as remarked in
\cite{GCT}, the ones we consider here apply without change to general
pseudo-Riemannian manifolds. The upshot of these remarks is that for every $x
\in U$, there is some $k(x) \geq 1$ so that the map:
\begin{align*}
	J^k_x : \Kill_0(M,x) &\rightarrow \Kill^k_0(M,x)  \\
		X &\mapsto j^k_x(X),
\end{align*}
is a linear isomorphism for every $k \geq k(x)$. Note that in this case, for
the usual brackets in $\Kill_0(M,x)$ and the brackets $[\cdot,\cdot]^k$ in
$\Kill^k_0(M,x)$ considered above, the map $J^k_x$ is a Lie algebra
anti-isomorphism.

For $S_0$ and $U$ as above, consider the dense subset $S = S_0 \cap U \subset
M$. Next choose $x \in S$ and $k \geq \max(k(x),2)$. Then, the map $J^k_x$
is a Lie algebra anti-isomorphism, and there exists a Lie algebra homomorphism
$\widehat{\rho}^k_x : \g \rightarrow \Kill^k_0(M,x)$ satisfying~(*).
If we let $\rho_x = - (J^k_x)^{-1}\circ \widehat{\rho}^k_x : \g \rightarrow
\Kill_0(M,x)$, then $\rho_x$ defines a Lie algebra homomorphism such
that:
$$
	j^{k-1}_x([\rho_x(X),Y^*]) = j^{k-1}_x([X,Y]^*),
$$
for every $X,Y \in \g$. For this, we have used (*) and the definition of
$\theta_x$ from Lemma~\ref{lemma-jets}. Since $k-1 \geq 1$ and because
germs of Killing fields are determined by any jet of order at least $1$, we
conclude that, at our chosen point $x$, $\rho_x$ satisfies (from our statement)
(1), (2) and the identity in (3) in a neighborhood of $x$ with $\widetilde{M}$
replaced with $M$.
The identity in (3) now proves that every element of $\g(x)$ preserves the
tangent bundle to $\OO$ in a neighborhood of $x$: i.e.~the corresponding Lie
derivatives map sections of $T\OO$ into sections of $T\OO$. By Proposition~2.2
of \cite{Molino} we conclude that the local flows of the elements of $\g(x)$
preserve $\OO$ as well in a neighborhood of $x$. Since the elements of $\g(x)$
are Killing fields, we conclude that they (and their local flows) also preserve
the normal bundle $T\OO^\perp$ in a neighborhood of $x$. This completes the
proof of our statement for the dense subset $S \subset M$ and for
$\widetilde{M}$ replaced with $M$ in (1)--(4). Finally, this yields the
statement for $\widetilde{M}$ for the dense subset which is the inverse image
of $S$ under the covering map since such map is a local isometry.
\end{proof}

\begin{remark}\label{remark-prop-g(x)}
The conclusions of Proposition~\ref{prop-g(x)} hold without change for some
dense subset $S \subset M$ by replacing $\widetilde{M}$ with $M$ in
(1)--(4). In fact, our proof first obtains the required Killing fields on $M$
which are then translated into corresponding ones on $\widetilde{M}$.
\end{remark}

In this work, we will be dealing with and interested in the case where the
$G$-orbits in $M$ are nondegenerate submanifolds with respect to the
pseudo-Riemannian metric. In this case, the $\widetilde{G}$-orbits on
$\widetilde{M}$ are nondegenerate as well and we have a direct sum decomposition
$T\widetilde{M} = T\OO \oplus T\OO^\perp$. When this holds, we can consider the
$\g$-valued $1$-form $\omega$ on $\widetilde{M}$ that is given, at every $x \in
\widetilde{M}$, by the composition $T_x\widetilde{M} \rightarrow
T_x\OO\cong\g$, where the first map is the natural projection and the second
map is the natural isomorphism. From this, we then define the $\g$-valued
$2$-form given by $\Omega = d\omega|_{\wedge^2 T\OO^\perp}$. The following
result will provide us with a criterion, in terms of $\Omega$, for the normal
bundle $T\OO^\perp$ to be integrable. At the same time, we relate $\Omega$ with
the $\g$-module structures from Proposition~\ref{prop-g(x)}. This result is
very well known and it is essentially contained in \cite{Gromov},
\cite{Hernandez} and \cite{Quiroga-Annals}, but we include its proof here for
the sake of completeness.

\begin{lemma}\label{lemma-omega-Omega}
Let $G$, $M$ and $S$ be as in Proposition~\ref{prop-g(x)}. If we assume that the
$G$-orbits are nondegenerate, then:
\begin{enumerate}
\item For every $x \in S$, the maps $\omega_x : T_x\widetilde{M} \rightarrow
\g$ and $\Omega_x : \wedge^2 T_x\OO^\perp \rightarrow \g$ are both homomorphisms
of $\g$-modules, for the $\g$-module structures from
Proposition~\ref{prop-g(x)}.

\item The normal bundle $T\OO^\perp$ is integrable if and only if $\Omega = 0$.
\end{enumerate}
\end{lemma}
\begin{proof}
In the rest of the proof, let $X \in \g $ and $x \in S$ be fixed but
arbitrarily given. 

For $Z$ a vector field over $\widetilde{M}$, let $Z^\top, Z^\perp$ be its $T\OO$
and $T\OO^\perp$ components, respectively. Since $\rho_x(X)$ is a Killing field
preserving $\OO$ and $T\OO^\perp$ it follows that:
\begin{align*}
	[\rho_x(X),Z]^\top &= [\rho_x(X),Z^\top], \\
	[\rho_x(X),Z]^\perp &= [\rho_x(X),Z^\perp].
\end{align*}
Denote with $\alpha : T_x\OO \rightarrow \g$ the inverse map of $X \mapsto
X^*_x$. Then we have:
\begin{align*}
	\omega_x(X\cdot Z_x) &= \omega_x([\rho_x(X),Z]_x)  \\
			&= \alpha([\rho_x(X),Z^\top]_x)  \\
			&= \alpha([\rho_x(X),\omega(Z)^*]_x)  \\
			&= \alpha([X,\omega(Z)]^*_x)  \\
			&= [X,\omega_x(Z)]  \\
			&= X\cdot \omega_x(Z_x),
\end{align*}
thus showing that $\omega_x$ is a homomorphism of $\g$-modules. Here we used in
the second and third identities the definition of $\omega$, and in the fourth
identity the formula from Proposition~\ref{prop-g(x)}(3); the rest follows from
the definition of the $\g$-module structures involved.

Next, observe that for every pair of sections $Z_1,Z_2$ of $T\OO^\perp$ we have:
\begin{align*}
	\Omega(Z_1\wedge Z_2) &= Z_1(\omega(Z_2)) - Z_2(\omega(Z_1)) 
			-\omega([Z_1,Z_2])  \\
			&= -\omega([Z_1,Z_2]),
\end{align*}
which clearly implies (2).

Now let $u,v \in T_x\OO^\perp$ be given and choose $U,V$ sections of
$T\OO^\perp$ extending them, respectively. Hence, using that $\omega$ is a
homomorphism of $\g$-modules, the Jacobi identity and the above expression
relating $\Omega$ and $\omega$, we obtain:
\begin{align*}
	\Omega_x(X\cdot(u\wedge v)) &= \Omega_x((X\cdot u)\wedge v) 
			+ \Omega_x(u \wedge (X\cdot v))  \\
		&= \Omega_x([\rho_x(X),U]\wedge V)
			+ \Omega_x(U\wedge [\rho_x(X),V])  \\
		&= -\omega_x([[\rho_x(X),U],V])
			- \omega_x([U,[\rho_x(X),V]])  \\
		&= -\omega_x([\rho_x(X),[U,V]])  \\
		&= -\omega_x(X\cdot [U,V]_x)  \\
		&= -[X,\omega_x([U,V])]  \\
		&= [X,\Omega_x(U\wedge V)] \\
		&= X \cdot \Omega_x(u\wedge v),
\end{align*}
thus showing that $\Omega_x$ is a homomorphism of $\g$-modules. Note that we
have used that both $[\rho_x(X),U], [\rho_x(X),V]$ are sections of $T\OO^\perp$.
\end{proof}

The following result allows us to relate the metric $T\OO$ coming from $M$
to suitable metrics on $G$. The proof presented here is due to Gromov
(see~\cite{Gromov}) and provides our first application of
Proposition~\ref{prop-g(x)}.

\begin{lemma}\label{lemma-metric-on-orbits}
Suppose that the $G$-action on $M$ has a dense orbit and preserves a finite
volume pseudo-Riemannian metric. Then, for every $x\in M$
and with respect to the natural isomorphism $\g \cong T_x\OO$, the metric of $M$
restricted to $T_x\OO$ defines an $\Ad(G)$-invariant symmetric bilinear form on
$\g$ independent of the point $x$.
\end{lemma}
\begin{proof}
With the above mentioned trivialization of $T\OO$, the metric $h$ on $M$
restricted to the orbits and pulled back to $M\times \g$ yields a map:
\begin{align*}
	\psi : M &\rightarrow \g^* \otimes \g^* \\
		x &\mapsto B_x
\end{align*}
where $B_x(X,Y) = h_x(X^*_x, Y^*_x)$.

By Remark~\ref{remark-prop-g(x)}, there is a dense subset $S \subset
M$ so that the conclusions of Proposition~\ref{prop-g(x)} are satisfied for
every $x \in S$ with the tangent spaces and Killing fields of $\widetilde{M}$
replaced by those of $M$. Hence, for every $x\in S$ the inner product $h_x$ is
preserved, in the sense of Proposition~\ref{prop-g(x)}(4), by the Killing
vector fields that belong to $\g(x)$. In particular, for every $x\in S$ and
$X,Y,Z \in \g$ we have:
$$
	h_x([\rho_x(X),Y^*]_x,Z^*_x) = - h_x(Y^*_x,[\rho_x(X),Z^*]_x),
$$
which, by Proposition~\ref{prop-g(x)}(3) yields:
$$
	h_x([X,Y]^*_x,Z^*_x) = - h_x(Y^*_x,[X,Z]^*_x).
$$
In other words, for every $x \in S$ we have:
$$
	B_x([X,Y],Z) = - B_x(Y,[X,Z]),
$$
for all $X,Y,Z \in \g$. This implies that $\psi(x) = B_x$ is an
$\Ad(G)$-invariant form on $\g$ for every $x \in S$. By the density of $S$ in
$M$, we conclude that the image of $\psi$ lies in the set of $\Ad(G)$-invariant
forms.

On the other hand, at every $x \in M$ and for $g\in G$, $X, Y \in \g$ we have:
\begin{align*}
	\psi(gx)(X,Y) &= h_{gx}(X^*_{gx},Y^*_{gx})  \\
		&= h_x(dg^{-1}_{gx}(X^*_{gx}),dg^{-1}_{gx}(Y^*_{gx}))  \\
		&= h_x(\Ad(g^{-1})(X)^*_x,\Ad(g^{-1})(Y)^*_x)  \\
		&= \psi(x)(\Ad(g^{-1})(X),\Ad(g^{-1})(Y)),
\end{align*}
which shows that $\psi$ is $G$-equivariant. Note that we used in the second
identity that $G$ preserves the metric, and in the third identity we used the
remarks at the beginning of this section.

The $G$-equivariance of $\psi$ and the fact that its image lies in $G$-fixed
points implies that $\psi$ is $G$-invariant. Then, the result follows from the
existence of a dense $G$-orbit.
\end{proof}

As a simple application of the previous result, we prove the nondegeneracy of
the orbits when the manifold acted upon has a suitably bounded dimension.

\begin{lemma}\label{lemma-nondeg}
Suppose that the $G$-action on $M$ has a dense orbit and preserves a finite
volume pseudo-Riemannian metric. If $\dim(M) < 2 \dim(G)$, then the $G$-orbits
are nondegenerate with respect to the metric on $M$.
\end{lemma}
\begin{proof}
By Lemma~\ref{lemma-metric-on-orbits}, for every $x \in M$ the metric
restricted to $T_x\OO$ corresponds to an $\Ad(G)$-invariant form in $\g$. The
kernel of such a form is an ideal and so the metric $h_x$ restricted to
$T_x\OO$ is either nondegenerate or zero.

Suppose that $h_x$ is zero when restricted to $T_x\OO$ for some $x \in M$. Then,
$T_x\OO$ lies in the null cone of $T_x M$ for the metric $h_x$. Hence,
for $(m_1,m_2)$ the signature of $M$, we have $\dim(G) =
\dim(T_x\OO) \leq \min(m_1,m_2)$. And this implies $2\dim(G) \leq m_1 +
m_2 = \dim(M)$, which is impossible.
\end{proof}

\section{Proof of Theorem~\ref{thm-int-normal} and some consequences}
\label{section-int-normal}
We start this section by proving Theorem~\ref{thm-int-normal}.

\begin{proof}[Proof of Theorem~\ref{thm-int-normal}]
Assuming that $T\mathcal{O}^\perp$ is integrable, let $\mathcal{F}$ be the
induced foliation. We will first prove that $\mathcal{F}$ is totally geodesic,
i.e.~its leaves are totally geodesic submanifolds of $M$. We will denote with
$h$ the metric on $M$ preserved by $G$.

First note that, if $Y,Z$ are local sections of $T\OO^\perp$ that preserve the
foliation, then we have for every $X \in \g$:
$$
	X^*(h(Y,Z)) = h([X^*,Y],Z) + h(Y,[X^*,Z]) = 0,
$$
because our choices imply that $[X^*,Y]$ and $[X^*,Z]$ are section of $T\OO$.
In particular, for every $Y,Z$ as above the function $h(Y,Z)$ is constant along
the $G$-orbits. In the notation of \cite{Molino}, we conclude that $h$ is a
bundle-like metric for the foliation $\OO$. Hence, by the remarks in page~79 of
\cite{Molino} it follows that $h$ induces a transverse metric to the foliation
$\OO$. By the construction of such transverse metric from $h$ and the arguments
in the proof of Proposition~3.2 of \cite{Molino}, it is easy to conclude that
the foliation $\OO$ is given by pseudo-Riemannian submersions that define the
transverse metric. More precisely, at every point in $M$ there is an open subset
$U$ of $M$ and a pseudo-Riemannian submersion $\pi : U \rightarrow B$ such that
the fibers of $\pi$ define the foliation $\mathcal{O}$ restricted to $U$. We
observe that the results of \cite{Molino} are stated for Riemannian metrics, but
those that we use here extend without change to arbitrary pseudo-Riemannian
metrics.

We will now use the properties of the structural equations from
\cite{ONeill-submersion} for a pseudo-Riemannian submersion $\pi : U \rightarrow
B$ as above. Again, the results in \cite{ONeill-submersion} are stated for
Riemannian submersions, but the ones that we will use are easily seen to hold
for pseudo-Riemannian submersions as well. For $\pi$ as above, let
$A$ be the associated fundamental tensor defined in \cite{ONeill-submersion}.
In particular, by the definition of $A$, the second fundamental form for the
leaves of $\mathcal{F}$ is given by $A_X Y$, for $X,Y$ tangent to
$\mathcal{F}$. But by Lemma~2 of \cite{ONeill-submersion} we have for $X,Y$
tangent to $\mathcal{F}$ the identity $A_X Y = \frac{1}{2}[X,Y]^\top$, where
$Z^\top$
denotes the projection of $Z$ onto $T\mathcal{O}$. Hence, the integrability of
$T\OO^\perp$ to $\mathcal{F}$ shows that $A$ vanishes on vector fields tangent
to $\mathcal{F}$, thus showing that the leaves of $\mathcal{F}$ are totally
geodesic.

Choose a leaf $N$ of $\mathcal{F}$. Then, one can prove fairly easy that every
geodesic in $M$ which is tangent at some point to $N$ remains in $N$ for every
value of the parameter of the geodesic; this uses the fact that $N$ is a maximal
integral submanifold of $T\OO^\perp$ and that the leaves of $\mathcal{F}$ are
totally geodesic. Hence, the completenes of $M$ implies that of $N$.

For our chosen leaf $N$ of $\mathcal{F}$, consider the $G$-action map
restricted to $G\times N$. This defines a smooth map $\varphi : G \times N
\rightarrow M$ which is $G$-equivariant. By Lemma~\ref{lemma-metric-on-orbits},
it follows easily that $\varphi$ is a local isometry for $G \times N$ endowed
with the product metric where $G$ carries a suitable bi-invariant metric. In
particular, $G\times N$ is complete and so we conclude from Corollary~29 in page
202 from \cite{ONeill-book} that $\varphi$ is an isometric covering map. Hence,
the universal covering map of $M$ is given by $\widetilde{\varphi} :
\widetilde{G} \times \widetilde{N} \rightarrow M$ and the
$\widetilde{G}$-action on $M$ lifted to $\widetilde{G} \times \widetilde{N}$ is
the left action on the first factor.

We now claim that $\pi_1(M) \subset \Iso(\widetilde{G}) \times
\Iso(\widetilde{N})$, i.e.~that every element in $\pi_1(M)$ preserves the
factors in the product $\widetilde{G} \times \widetilde{N}$. To see this, let
$\gamma \in \pi_1(M)$ be given with $\gamma = (\gamma_1,\gamma_2)$ its
component decomposition. Observe that in $\widetilde{G} \times \widetilde{N}$
we have:
$$
	T_{(g,x)}\OO = T_g\widetilde{G}, \quad 
		T_{(g,x)}\OO^\perp = T_x \widetilde{N},
$$
for every $(g,x) \in \widetilde{G} \times \widetilde{N}$. Since $\gamma$
commutes with the $\widetilde{G}$-action, it preserves both $T\OO$ and
$T\OO^\perp$ and so:
\begin{align*}
	d\gamma(u) &= d\gamma_1(u) + d\gamma_2(u) \in T\OO  \\
	d\gamma(v) &= d\gamma_1(v) + d\gamma_2(v) \in T\OO^\perp, 
\end{align*}
for every $u\in T\widetilde{G}$ and $v \in T\widetilde{N}$. We conclude that
$d\gamma_2(T\OO) = 0$ and $d\gamma_1(T\OO^\perp) = 0$, which implies that
$\gamma_1$ is independent of $\widetilde{N}$ and $\gamma_2$ is independent of
$\widetilde{G}$. This yields our claim about $\pi_1(M)$.

On the other hand, since $\widetilde{G}$ carries a bi-invariant metric, by the
results from Section~4 of \cite{Quiroga-Annals} we know that the connected
component of the identity of $\Iso(\widetilde{G})$ is given by
$\Iso_0(\widetilde{G}) = L(\widetilde{G})R(\widetilde{G})$ (the left and right
translations) and that it is a finite index subgroup of $\Iso(\widetilde{G})$.
Hence, the group $\Lambda = \pi_1(M) \cap (\Iso_0(\widetilde{G}) \times
\Iso(\widetilde{N}))$ has finite index in $\pi_1(M)$, and so the induced map
$(\widetilde{G} \times \widetilde{N})/ \Lambda \rightarrow M$ is a finite
covering. Moreover, every $\gamma \in \Lambda$ can be written as $\gamma =
(L_{g_1} R_{g_2}, \gamma_2)$, where $g_1, g_2 \in \widetilde{G}$ and $\gamma_2
\in \Iso(\widetilde{N})$, and since such $\gamma$ commutes with the
$\widetilde{G}$-action we conclude that:
$$
	(g_1gg_2,\gamma_2(x)) = (L_{g_1} R_{g_2},\gamma_2)(g(e,x)) =
    		g((L_{g_1} R_{g_2},\gamma_2)(e,x)) = (gg_1g_2,\gamma_2(x))
$$
for every $g\in \widetilde{G}$ and $x\in \widetilde{N}$, which implies $g_1
\in Z(\widetilde{G})$. Hence, $L_{g_1} = R_{g_1}$ and then $\gamma \in
R(\widetilde{G})\times \Iso(\widetilde{N})$, thus showing that $\Lambda \subset
R(\widetilde{G})\times \Iso(\widetilde{N})$.

Also note that the covering map $G \times \widetilde{N} \rightarrow
M$ realizes $\pi_1(G) \subset \Lambda$, which induces a covering map $G\times
\widetilde{N} \rightarrow (\widetilde{G} \times \widetilde{N})/\Lambda$. We
claim that $G\times \widetilde{N} \rightarrow (\widetilde{G} \times
\widetilde{N})/\Lambda$ is a normal covering map. For this
we need to check that $\pi_1(G)$ is a normal subgroup of $\Lambda$
under the inclusion $z \mapsto (R_z,e)$. But by the above remarks, every $\gamma
\in \Lambda$ can be written as $\gamma = (R_{g_1}, \gamma_2)$, where
$g_1 \in \widetilde{G}$ and $\gamma_2 \in \Iso(\widetilde{N})$, from which
we obtain:
$$
    \gamma (R_z,e) \gamma^{-1}
        = (R_{g_1} R_z R_{g^{-1}_1},\gamma_2\gamma^{-1}_2)
        = (R_z,e)
$$
since $\pi_1(G)$ is central in $\widetilde{G}$. It follows that the group of
deck transformations for $G\times \widetilde{N} \rightarrow (\widetilde{G}
\times \widetilde{N})/\Lambda$ is given by the group $\Gamma = \Lambda /
\pi_1(G)$ and that we also have $(G \times \widetilde{N})/\Gamma =
(\widetilde{G} \times \widetilde{N})/\Lambda$.

From the above, we conclude that $\Gamma \subset R(G) \times
\Iso(\widetilde{N}) = G \times \Iso(\widetilde{N})$ is a group of deck
transformations of $G\times \widetilde{N} \rightarrow M$ that induces a
$G$-equivariant finite covering map $(G\times \widetilde{N})/\Gamma \rightarrow
M$ that satisfies the properties required to obtain
Theorem~\ref{thm-int-normal}.
\end{proof}

The integrability of the normal bundle can be ensured for suitable relations
between the Lie group $G$ and the geometry of the manifold on which it acts,
thus providing the following results. In what follows, the signature of $G$, as
a pseudo-Riemannian manifold, is always considered with respect to a
bi-invariant metric. Note that if $(n_1,n_2)$ is the signature of some
bi-invariant metric on $G$, then the signature of any other bi-invariant metric
is either $(n_1,n_2)$ or $(n_2,n_1)$.

\begin{corollary}\label{cor-int-normal-rank}
Suppose that the $G$-action on $M$ has a dense orbit and preserves a finite
volume complete pseudo-Riemannian metric. If the $G$-orbits are nondegenerate
and either one of the following holds:
\begin{enumerate}
\item there is no Lie algebra embedding of $\g$ into $\so(T_x\OO^\perp)$ for
every $x \in M$, or

\item for $n_0 = \min(n_1,n_2)$ and $m_0 = \min(m_1,m_2)$, where $(n_1,n_2)$
and $(m_1,m_2)$ are the signatures of $G$ and $M$, respectively, we have
$\rank_\R(\g) > m_0 - n_0$,
\end{enumerate}
then the conclusion of Theorem~\ref{thm-int-normal} holds.
\end{corollary}
\begin{proof}
Let us consider a subset $S \subset \widetilde{M}$ given as in
Proposition~\ref{prop-g(x)}.

First suppose that condition (1) is satisfied. Then, for every $x \in S$, the
$\g$-module structure on $T_x\OO^\perp$ given by Proposition~\ref{prop-g(x)}(4)
is trivial. By Lemma~\ref{lemma-omega-Omega}(1), being a homomorphism of
$\g$-modules, the map $\Omega_x$ is trivial for every $x \in S$. Since $S$ is
dense, we conclude that $\Omega = 0$ and so $T\OO^\perp$ is integrable by
Lemma~\ref{lemma-omega-Omega}(2). Hence, Theorem~\ref{thm-int-normal} can be
applied.

Let us now assume that (2) holds. Note that by
Lemma~\ref{lemma-metric-on-orbits}, the signature of $T\OO$ is either
$(n_1,n_2)$ or $(n_2,n_1)$. If we let $(k_1,k_2)$ be the signature of
$T\OO^\perp$, then it is easily seen that:
$$
	\min(k_1,k_2) \leq m_0 - n_0.
$$
Since the real rank of $\so(T_x\OO^\perp)$ is precisely $\min(k_1,k_2)$ the
result in this case follows from the first part.
\end{proof}

\begin{corollary}\label{cor-int-normal-reps}
Suppose that the $G$-action on $M$ has a dense orbit and preserves a finite
volume complete pseudo-Riemannian metric. Let $n$ be the dimension of the
smallest $\g$-module $V$ such that $\wedge^2 V$ contains a $\g$-submodule
isomorphic to $\g$. If $\dim(M) < \dim(G) + n$, then the conclusion of
Theorem~\ref{thm-int-normal} holds.
\end{corollary}
\begin{proof}
First observe that the Lie brackets in $\g$ define a surjective homomorphism of
$\g$-modules from $\wedge^2 \g$ onto $\g$. Hence, $\wedge^2 \g$ contains a
submodule isomorphic to $\g$, thus implying that $n \leq \dim(\g)$. Then,
Lemma~\ref{lemma-nondeg} shows that the $G$-orbits in $M$ are nondegenerate
with respect to the metric of $M$.

Let $S\subset \widetilde{M}$ be a subset given as in
Proposition~\ref{prop-g(x)}. From our hypotheses and since the $G$-orbits are 
nondegenerate, we have $\dim(T_x\OO^\perp) < n$. In particular, for every $x \in
S$ and for the $\g$-module structure from Proposition~\ref{prop-g(x)}(4),
$\wedge^2 T_x\OO^\perp$ does not contain a $\g$-submodule isomorphic to $\g$.
Hence, Lemma~\ref{lemma-omega-Omega}(1) implies that $\Omega_x = 0$
for every $x \in S$ and so that $\Omega = 0$. By
Lemma~\ref{lemma-omega-Omega}(2) the bundle $T\OO^\perp$ is integrable and so
Theorem~\ref{thm-int-normal} can be applied.
\end{proof}

On a compact manifold with Riemannian normal bundle, we can obtain the
conclusions of Theorem~\ref{thm-int-normal} without having to a priori assume
that the manifold is complete and that the normal bundle is integrable.

\begin{corollary}\label{cor-int-compact-Riemannian}
Suppose that the $G$-action on $M$ has a dense orbit and preserves a
pseudo-Riemannian metric. If $M$ is compact, the $G$-orbits are nondegenerate
and the normal bundle $T\OO^\perp$ is Riemannian, then the conclusions of
Theorem~\ref{thm-int-normal} hold. Moreover, we can assume that $\widetilde{N}$
is Riemannian homogeneous and that $\Gamma \subset G \times
\Iso_0(\widetilde{N})$.
\end{corollary}
\begin{proof}
The proof is a refinement of that of Theorem~\ref{thm-int-normal}, so we will
follow the notation of the latter.

First observe that the integrability of $T\OO^\perp$ follows from the proof of
Corollary~\ref{cor-int-normal-rank}(1) since $\so(T_x\OO^\perp)$ is compact for
every $x \in \widetilde{M}$. In particular, we have the hypotheses
of Theorem~\ref{thm-int-normal} except for the completeness of $M$. With this
respect, it is easy to check that the compactness of $M$ and the fact that
$T\OO^\perp$ is Riemannian imply that the geodesics in $M$ perpendicular to the
$G$-orbits are complete. This completeness is enough for the rest of the
arguments in the proof of Theorem~\ref{thm-int-normal} to apply.

Finally, we observe that the existence of a dense $G$-orbit implies that
$\widetilde{N}$ has a dense orbit by its local isometries and, being Riemannian,
we conclude that it is homogeneous. The latter follows from the infinitesimal
characterization of homogeneous Riemannian manifolds obtained in \cite{Singer},
and the fact that the orthogonal group is compact for definite metrics; we
refer to \cite{Quiroga-Annals} for further details. But for a
homogeneous Riemannian manifold the group of isometries has finitely many
connected components, and so we can intersect $\Gamma$ with $G \times
\Iso_0(\widetilde{N})$ to obtain the last claim after passing to
a finite covering.
\end{proof}

\section{Manifolds with a transverse Riemannian structure: Proof of
Theorem~\ref{theorem-trans-Riem-equiv}} \label{section-Riem}
In this section we will characterize those actions that preserve a Riemannian
structure transverse to the $G$-orbits. We start with the following basic
result relating transverse geometric structures for the foliation by $G$-orbits
with geometric structures on the normal bundle to such orbits.

\begin{lemma}\label{lemma-transverse-geometric}
Suppose that the $G$-action on $M$ has a dense orbit and preserves a finite
volume pseudo-Riemannian metric. Also assume that the $G$-orbits are
nondegenerate. If $H$ is a subgroup of $\mathrm{GL}(k,\R)$, where $k = \dim(M) -
\dim(G)$, then there is a one-to-one correspondence between the $G$-invariant
$H$-reductions of $L(TM/T\OO)$ and the $G$-invariant $H$-reductions of
$L(T\OO^\perp)$. In particular, every transverse $H$-structure for the foliation
$\OO$ by $G$-orbits induces a $G$-invariant $H$-reduction of $L(T\OO^\perp)$.
\end{lemma}
\begin{proof}
From the decomposition $TM = T\OO\oplus T\OO^\perp$ we obtain a natural
$G$-equivariant isomorphism $TM/T\OO \rightarrow T\OO^\perp$ which clearly
yields the first claim.

Next, we recall that a transverse $H$-structure to the foliation $\OO$ is given
by a reduction $P$ of $L(TM/T\OO)$ which is invariant under the local flows of
vectors fields tangent to the foliation $\OO$. In particular, $P$ is invariant
under the $G$-action on $L(TM/T\OO)$ thus showing the last claim.
\end{proof}

The following result is obtained by applying the main theorems
from~\cite{Quiroga-Annals}. 

\begin{proposition}\label{prop-from-Annals}
Suppose that the $G$-action on $M$ has a dense orbit and preserves a
pseudo-Riemannian metric. Also assume that $M$ is compact. If the normal bundle
to the orbits $T\OO^\perp$ is Riemannian, then there exist:
\begin{enumerate}
\item a finite covering map $\widehat{M} \rightarrow M$,

\item a connected Lie group $H$ with a compact subgroup $K$, and

\item a discrete cocompact subgroup $\Gamma \subset G\times H$ such that
$G\Gamma$ is dense in $G\times H$, 
\end{enumerate}
for which the $G$-action on $M$ lifts to $\widehat{M}$ so that $\widehat{M}$ is
$G$-equivariantly diffeomorphic to $(G \times K\backslash
H)/\Gamma$. Furthermore, if $G$ has finite center and real rank at least $2$,
then we can assume that $G\times H$ is a finite center isotypic
semisimple Lie group and $\Gamma$ is an irreducible lattice.
\end{proposition}
\begin{proof}
Let us denote with $(m_1,m_2)$ and $(n_1,n_2)$ the signatures of $M$ and $G$,
respectively. Also, let us denote $m_0 = \min(m_1,m_2)$ and $n_0 =
\min(n_1,n_2)$. Observe that since $T\OO^\perp$ is Riemannian, the bundle $T\OO$
is nondegenerate. Hence, by Lemma~\ref{lemma-metric-on-orbits} the signature of
$T\OO$ is either $(n_1,n_2)$ or $(n_2,n_1)$; this is because the signature of
every $\Ad(G)$-invariant metric on $\g$ is the signature of either the Killing
form or its negative. Without loss of generality, we can assume that the
signature of $T\OO$ is precisely $(n_1,n_2)$.

Theorems~A and B from \cite{Quiroga-Annals} provide precisely our required
conclusions when $n_0 = m_0$ holds, except for explicitly stating the
property that $G\Gamma$ is dense in $G\times H$. The latter is obtained as
follows. The conclusion of Theorem~A from \cite{Quiroga-Annals} ensures the
existence of a $G$-invariant ergodic smooth measure, which in turn implies the
existence of a dense $G$-orbit in $(G\times H)/\Gamma$. From this, it is easily
seen that we necessarily have the density of the $G$-orbit of $(e,e)\Gamma$ in
$(G\times H)/\Gamma$. Since the map $G \times H \rightarrow (G\times H)/\Gamma$
is a covering, we conclude that the inverse image under this covering of such
orbit, which is $G\Gamma$, is dense in $G \times H$.

We now consider the possibilities for $n_0$. We will now assume that the
signatures are given in the form $(+,-)$. If $n_0 = n_2$, then the fact that
$T\OO^\perp$ is Riemannian implies that $(m_1,m_2) = (n_1+k,n_2)$, where $k$ is
the rank of $T\OO^\perp$. In particular, $n_0 = m_0$ in this case and the
result follows from \cite{Quiroga-Annals}.

In case $n_0 = n_1$, we can replace the metric on $T\OO^\perp$ by its negative
to obtain a new $G$-invariant metric for which the signature of $M$ is
$(m_1,m_2) = (n_1,n_2+k)$, where $k$ is again the rank of $T\OO^\perp$. Hence,
for this new metric $n_0 = m_0$ and so the result follows also from
\cite{Quiroga-Annals}.
\end{proof}

We now proceed to prove Theorem~\ref{theorem-trans-Riem-equiv}.

\begin{proof}[Proof of Theorem~\ref{theorem-trans-Riem-equiv}]
Assume that (1) holds. Let us take $\widehat{M} = (G \times K\backslash
H)/\Gamma$ and verify that it satisfies (2). 

Note that, with respect to the quotient map $G \times H \rightarrow (G \times
H)/\Gamma$, the set $G\Gamma$ projects onto the $G$-orbit of the class of
the identity. Hence, there is a dense $G$-orbit in $\widehat{M}$.

Next, endow $G \times H$ with the product metric given by a bi-invariant metric
on $G$ and a Riemannian metric on $H$ which is left $K$-invariant and right
$H$-invariant. The latter exists because $K$ is compact. This induces a
pseudo-Riemannian metric on $G \times K \backslash H$ which is left
$G$-invariant and right $\Gamma$-invariant. Note that the $G$-orbits are
nondegenerate with normal bundle given by the tangent bundle to the factor
$K\backslash H$. Furthermore, by construction the projection $G \times
K\backslash H \rightarrow K\backslash H$ is a pseudo-Riemannian submersion and
so it defines a transverse Riemannian structure for the foliation by
$G$-orbits. Since the left $G$-action and the $\Gamma$-action commute with each
other, the transverse Riemannian structure on $G \times K \backslash H$ induces
a corresponding one on the double coset $(G \times K \backslash H)/\Gamma$
which is $G$-invariant. This provides the geometric structures required by (2).

Let us now assume that (2) holds, so that in particular $\widehat{M}$ has
a dense $G$-orbit and carries the indicated geometric structures. Since the
$G$-orbits in $\widehat{M}$ are nondegenerate then we have an orthogonal
decomposition $T\widehat{M} = T\OO \oplus T\OO^\perp$ which is invariant under
$G$. The existence of a transverse Riemannian structure for the foliation by
$G$-orbits in $\widehat{M}$ and Lemma~\ref{lemma-transverse-geometric} imply the
existence of a $G$-invariant Riemannian metric on $T\OO^\perp$. Hence, if we
replace with such Riemannian metric the metric on $T\OO^\perp$ induced from
$\widehat{M}$, we obtain a new $G$-invariant pseudo-Riemannian metric on
$\widehat{M}$ for which $T\OO^\perp$ is Riemannian. Hence, the $G$-action on
$\widehat{M}$ satisfies the hypotheses of Proposition~\ref{prop-from-Annals},
and the latter provides a finite covering of $\widehat{M}$, and thus of $M$,
with the properties required by (1).

Finally, let us assume that $G$ has finite center and real rank at least $2$.
First observe that (1) replaced by the conditions in the last part of the 
statement still implies (2). The only nontrivial property to check is the
existence of a dense $G$-orbit in $(G \times K \backslash H)/\Gamma$ for
$\Gamma$ an irreducible lattice. But in this case, $G\Gamma$ is dense in
$G\times H$ (see for example Lemma~6.3 from \cite{Quiroga-Annals}) which yields
the required dense orbit. On the other hand, that (2) implies (1) with the
properties from the last part of the statement is a consequence of
Proposition~\ref{prop-from-Annals}.
\end{proof}

\section{Actions on Lorentzian manifolds}
\label{section-Lorentz}
In this section, we present a characterization of the $G$-actions on $M$
preserving a Lorentzian metric.

\begin{theorem}\label{theorem-Lorentzian}
Let $G$ be a connected noncompact simple Lie group acting faithfully on a
compact manifold $M$. Then the following conditions are equivalent:
\begin{enumerate}
\item The group $G$ is locally isomorphic to $\SL(2,\R)$ and there is a
$G$-equivariant finite covering map $(G \times K\backslash H)/\Gamma \rightarrow
M$ where $H$ is a connected Lie group with a compact subgroup $K$ and $\Gamma
\subset G \times H$ is a discrete cocompact subgroup such that $G\Gamma$ is
dense in $G \times H$.

\item There is a finite covering map $\widehat{M} \rightarrow M$ for which the
$G$-action on $M$ lifts to a $G$-action on $\widehat{M}$ with a dense orbit and
preserving a Lorentzian metric.
\end{enumerate}
\end{theorem}
\begin{proof}
Let us assume (1) and consider the finite covering $\widehat{M} = (G \times
K\backslash H)/\Gamma \rightarrow M$. Endow $G \times H$ with the product metric
given by the bi-invariant metric coming from the Killing form in $\g$ and a left
$K$-invariant and right $H$-invariant Riemannian metric on $H$. 
Then consider the induced
$G$-invariant metric on $\widehat{M}$. Since $G$ is locally isomorphic to
$\SL(2,\R)$, such metric on $\widehat{M}$ is Lorentzian. As in the proof of
Theorem~\ref{theorem-trans-Riem-equiv}, the density of $G\Gamma$ implies the
existence of a dense $G$-orbit in $\widehat{M}$. This proves (2).

Let us now assume that (2) holds. By Lemma~\ref{lemma-metric-on-orbits}, for
every $x \in M$, the metric in $T_x\OO$ corresponds to an $\Ad(G)$-invariant
symmetric bilinear form on $\g$. Such form is either $0$ or nondegenerate. The
former case implies the existence of a null tangent subspace of dimension at
least $3$, which is impossible; in particular, the $G$-orbits are nondegenerate.
Since $G$ is noncompact and simple, we conclude the existence of an
$\Ad(G)$-invariant form on $\g$ which is Lorentzian. But $\slinear(2,\R)$ is
the only simple Lie algebra admitting such a form, which implies that $G$ is
locally isomorphic to $\SL(2,\R)$. We also conclude that $T\OO^\perp$ is
Riemannian and so the rest of the claims in (1) follow from
Proposition~\ref{prop-from-Annals}.
\end{proof}

\section{Manifolds with a transverse parallelism: Proof of
Theorem~\ref{theorem-trans-parallelism-equiv}}\label{section-parallel}
The following result is an immediate consequence of the properties of Lie
foliations. A proof can be found in \cite{Quiroga-Annals}.

\begin{lemma}\label{lemma-fol-Lie-dense}
Let $X$ be a compact manifold carrying a foliation with a transverse Lie
structure. If the foliation has a dense leaf, then the lifted foliation to any
finite covering space of $X$ has a dense leaf as well.
\end{lemma}

We now obtain the next result which describes actions preserving a metric
and a transverse parallelism. Its proof is based on some of the arguments found
in \cite{Quiroga-AJM}.

\begin{proposition}\label{prop-from-AJM}
Suppose that the $G$-action on $M$ has a dense orbit and preserves a
pseudo-Riemannian metric. Also assume that $M$ is compact. If the foliation by
$G$-orbits is nondegenerate (with respect to the pseudo-Riemannian metric) and
carries a transverse parallelism, then there exist:
\begin{enumerate}
\item a finite covering map $\widehat{M} \rightarrow M$,

\item a connected Lie group $H$, and

\item a discrete cocompact subgroup $\Gamma \subset G\times H$ such that
$G\Gamma$ is dense in $G \times H$,
\end{enumerate}
for which the $G$-action on $M$ lifts to $\widehat{M}$ so that $\widehat{M}$ is
$G$-equivariantly diffeomorphic to $(G \times H)/\Gamma$. Furthermore, if $G$
has finite center and real rank at least $2$, then we can assume that $G\times
H$ is a finite center isotypic semisimple Lie group and $\Gamma$ is an
irreducible lattice.
\end{proposition}
\begin{proof}
By Lemma~\ref{lemma-transverse-geometric}, the transverse parallelism to the
$G$-orbits yields a $G$-invariant trivialization of $L(T\OO^\perp)$. Hence,
there is a family of $G$-invariant sections of $T\OO^\perp$, say $X_1, \dots,
X_k$, that defines a basis of $T\OO^\perp$ on every fiber. Let us consider the
$G$-invariant Riemannian metric on $T\OO^\perp$ for which these vector fields
are orthonormal at every point. Because of the orthogonal decomposition $TM = T
\OO \oplus T\OO^\perp$, if we replace the metric on $T\OO^\perp$ induced from
$M$ with the $G$-invariant Riemannian metric thus defined from the parallelism,
then we obtain a pseudo-Riemannian metric $h$ on $M$ which is $G$-invariant,
defines the same orthogonal complement $T\OO^\perp$ to the orbits and such that
this orthogonal complement is Riemannian. In the rest of this proof we will
consider $M$ endowed with this new metric $h$. Hence, by
Corollary~\ref{cor-int-compact-Riemannian} there is a simply connected
homogeneous Riemannian manifold $\widetilde{N}$ and a discrete cocompact
subgroup $\Gamma \subset G \times \Iso_0(\widetilde{N})$ for which there is a
$G$-equivariant finite covering map:
$$
	\widehat{M} = (G \times \widetilde{N})/\Gamma \rightarrow M.
$$

Note that by the proof of Corollary~\ref{cor-int-compact-Riemannian},
the normal bundle $T\OO^\perp$ is integrable. In particular, the vector fields
$X_i$ given above satisfy:
$$
	[X_i,X_j] = \sum_{r=1}^k f^r_{ij} X_r
$$
for some smooth functions $f^r_{ij}$ defined on $M$. The $G$-invariance of the
parallelism then implies that the functions $f^r_{ij}$ are $G$-invariant and so
constant because of the existence of a dense $G$-orbit. We conclude that the
linear span over $\R$ of the vector fields $X_1, \dots, X_k$ is a Lie algebra.

On the other hand, the fact that the vector fields $X_i$ are preserved by the
$G$-action is easily seen to imply that such fields are foliate in the notation
of \cite{Molino}. Hence, the parallelism $X_1,\dots, X_k$ defines a transverse
Lie structure for the foliation by $G$-orbits in $M$. 
Clearly, this transverse Lie structure induces a
corresponding one on $\widehat{M}$. Let $H$ be a simply connected Lie group that
models the transverse Lie structure on $\widehat{M}$ and consider a
corresponding development $D : \widetilde{G} \times \widetilde{N} \rightarrow
H$. In particular, $D$ is a submersion whose fibers have connected components
given precisely by subsets of the form $\widetilde{G} \times \{x\}$. 
Hence, we conclude that $D_{\{e\}\times \widetilde{N}} : \widetilde{N}
\rightarrow H$ defines a local diffeomorphism.

By the proofs of Corollary~\ref{cor-int-compact-Riemannian} and
Theorem~\ref{thm-int-normal}, the manifold $\widetilde{N}$ is the (isometric)
universal covering space of a leaf $N$ for the foliation in $M$ defined by
$T\OO^\perp$. Hence, if we let $\widetilde{X}_1, \dots, \widetilde{X}_k$ be the
pull backs to $\widetilde{N}$ of the restrictions of the fields $X_1, \dots,
X_k$ to $N$, then the induced fields on $\widetilde{G} \times \widetilde{N}$
(which we will denote with the same symbols) define the transverse Lie
structure on $\widetilde{G} \times \widetilde{N}$. Furthermore, from the
previous construction of the metric $h$, the metric on $\widetilde{N}$ is given
by the condition of the fields $\widetilde{X}_1, \dots, \widetilde{X}_k$ being
orthonormal at every point.
In particular, if we let $\h$ be the Lie algebra of $H$, then there is a
basis $v_1, \dots, v_k$ of $\h$ such that the development $D$ maps the vector
field $\widetilde{X}_i$ into $v_i$, for every $i=1,\dots,k$. Without loss of
generality we will assume that the transverse $H$-structures are modeled by
taking $H$ with its right translations, and so transverse Lie parallelisms are
modeled on right invariant vector fields. With this convention, each $v_i$ is
considered as a right invariant vector field on $H$. 

From the above remarks, if we endow $H$ with the right invariant Riemannian
metric for which the vector fields $v_1, \dots, v_k$ are orthonormal at every
point, then $D_{\{e\}\times \widetilde{N}} : \widetilde{N} \rightarrow H$ is a
local isometry. By Corollary~29 in page 202 of
\cite{ONeill-book} and since $\widetilde{N}$ is complete, we conclude that 
$D_{\{e\}\times \widetilde{N}} : \widetilde{N} \rightarrow H$ is an isometry
and so it induces in $\widetilde{N}$ a Lie group structure with respect to
which it is an isomorphism. Hence, we can replace the Riemannian manifold
$\widetilde{N}$ with $H$ carrying the above right invariant Riemannian
metric and assume that the natural projection $\widetilde{G} \times H
\rightarrow H$ is a development for the transverse Lie structure on
$\widehat{M}$. Moreover, by the definition of $\widehat{M}$, the space $G
\times \widetilde{N}$ covers $\widehat{M}$, and so we can assume that the
natural projection $\pi : G \times H \rightarrow H$ is also a development for
the
transverse Lie structure on $\widehat{M}$ with a corresponding holonomy
representation $\widehat{\rho} : \Gamma \rightarrow H$.

Since $\Gamma \subset G \times \Iso_0(\widetilde{N}) = G \times
\Iso_0(H)$, if we choose $\gamma \in \Gamma$, then we can write $\gamma =
(R_{\gamma_1},\gamma_2)$, where $\gamma_1 \in G$ and $\gamma_2 \in \Iso_0(H)$.
And so, the $\widehat{\rho}$-equivariance of $\pi$ yields for every $(g,x) \in
G\times H$:
$$
	\gamma_2(x) = \pi(g\gamma_1,\gamma_2(x)) 
		= \pi((g,x)\gamma)
		= \pi(g,x)\widehat{\rho}(\gamma)
		= xy
$$
where $y \in H$. It follows that $\gamma_2$ is given by a right translation by
an element in $H$, and so we have $\Gamma \subset G \times R(H) = G \times H$.
Since $\widehat{M} = (G\times H)/\Gamma \rightarrow M$ is a finite covering, by
Lemma~\ref{lemma-fol-Lie-dense} we conclude that $\widehat{M}$ has a dense
$G$-orbit and so $G\Gamma$ is dense in $G \times H$ as in the proof of
Proposition~\ref{prop-from-Annals}. This concludes the proof of the first part
of the statement.

For the last claim, if $G$ has finite center and real rank at least $2$, then
$H$ is semisimple by the main results from \cite{Zimmer-Lie}. Also, by modding
out by a suitable central group we can assume that $H$ has finite center (see
\cite{Quiroga-Annals}). The arguments at the end of Section~6 from
\cite{Quiroga-Annals} also prove that $\Gamma$ is an irreducible lattice.
Finally, recall that an irreducible lattice in a semisimple Lie group can only
exist if the group is isotypic.
\end{proof}

We can now prove our characterization of actions with a transverse parallelism.

\begin{proof}[Proof of Theorem~\ref{theorem-trans-parallelism-equiv}]
First, let us assume that (1) holds and consider the finite covering
$\widehat{M} = (G \times H)/\Gamma \rightarrow M$. As in the proof of
Theorem~\ref{theorem-trans-Riem-equiv} we conclude that $\widehat{M}$ has a
dense $G$-orbit. Note that the right invariant vector fields on $H$ define both
a transverse Lie structure on $G\times H$ for the foliation given by the factor
$G$, and a right $H$-invariant Riemannian metric on $H$. If we consider
the product pseudo-Riemannian metric on $G \times H$ using a bi-invariant
metric on $G$, then we also obtain a pseudo-Riemannian metric for which the
foliation given by the factor $G$ is nondegenerate. Both of these geometric
structures on $G \times H$ are left $G$-invariant and right $\Gamma$-invariant
and so descend to corresponding $G$-invariant geometric structures on
$\widehat{M}$ thus establishing (2).

If we now assume (2), then Proposition~\ref{prop-from-AJM} applied to
$\widehat{M}$ yields (1). The last claim is also a consequence of
Proposition~\ref{prop-from-AJM}.
\end{proof}

\section{Proof of Theorem~\ref{theorem-alghull-TOperp-irred}}
\label{section-alghull}
Let us consider $G$ and $M$ as in Proposition~\ref{prop-g(x)} with $S\subset
M$ a dense subset provided by this result and Remark~\ref{remark-prop-g(x)}.
With the notation of
Proposition~\ref{prop-g(x)}, the representation $\lambda_x\circ \rho_x$ leaves
invariant the subspace $T_x\OO^\perp$ thus defining its $\g$-module structure.
This induces a homomorphism $\g \rightarrow \so(T_x\OO^\perp)$ obtained by
restricting $\lambda_x\circ \rho_x$ to the submodule $T_x\OO^\perp$. By
composition with such homomorphism, we can consider $\Omega_x$ obtained from
Lemma~\ref{lemma-omega-Omega} as a map $\wedge^2 T_x\OO^\perp \rightarrow
\so(T_x\OO^\perp)$. 

For a vector space $W$ with inner product $\left<\cdot,\cdot\right>$, we will
say that a bilinear map $T : W\times W \rightarrow \mathfrak{gl}(W)$ is of
curvature type if it satisfies the following conditions for every $x,y,z,v,w \in
W$:
\begin{enumerate}
\item $T(x,y) = - T(x,y)$,

\item $\left<T(x,y)v,w\right> = - \left<v,T(x,y)w\right>$,

\item $T(x,y)z + T(y,z)x + T(z,x)y = 0$,

\item $\left<T(x,y)v,w\right> = \left<T(v,w)x,y\right>$.
\end{enumerate}
Note that (1) and (2) together are equivalent to $T$ inducing a map $\wedge^2 W
\rightarrow \so(W)$. Also, by the proof of Proposition~36 in page~75 of
\cite{ONeill-book} we know that (1), (2) and (3) together imply (4).

With the above notation, we have the following result.

\begin{lemma}\label{lemma-Omega-curvature}
Let $G$, $M$ and $S \subset M$ be as in Proposition~\ref{prop-g(x)} and
Remark~\ref{remark-prop-g(x)}.
For every $x \in S$, define the bilinear operation $[\cdot,\cdot]_0$ in
$T_x M$ by the assignments:
\begin{itemize}
\item $[X^*_x,Y^*_x]_0 = [X,Y]^*_x$, for every $X^*_x, Y^*_x \in T_x\OO$, $(X,Y
\in \g)$,

\item $[X^*_x,v]_0 = - [v,X^*_x]_0 = X(v)$ for every $X^*_x \in T_x\OO, v \in
T_x\OO^\perp$, $(X \in \g)$,

\item $[v_1,v_2]_0 = \Omega_x(v_1,v_2)^*_x$ for every $v_1, v_2 \in
T_x\OO^\perp$.
\end{itemize}
Then, $[\cdot,\cdot]_0$ yields a Lie algebra structure on $T_x M$
if and only if $\Omega_x$ is of curvature type when considered as a map 
$\wedge^2 T_x\OO^\perp \rightarrow \so(T_x\OO^\perp)$. In this case, the 
representation of $\g$ in $T_x M$, given by
Proposition~\ref{prop-g(x)}(4), preserves the Lie algebra structure of
$T_x M$.
\end{lemma}
\begin{proof}
Note that by the above remarks, $\Omega_x$ is of curvature type if and only
if it satisfies the above condition (3). Also observe that $[\cdot,\cdot]_0$
always defines a skew-symmetric bilinear form in $T_x M$. Hence, for
the first claim, we need to show that $[\cdot,\cdot]_0$ satisfies the Jacobi
identity if and only if $\Omega_x$ satisfies (3).

For the Jacobi identity to hold we only need to verify the following cases.

$\bullet$ $u_1, u_2, u_3 \in T_x\OO$. Note that $[\cdot,\cdot]_0$ maps
$T_x\OO\times T_x\OO$ into $T_x\OO$ in such a way that it defines a
skew-symmetric operation that corresponds to the Lie brackets of $\g$ under the
natural isomorphism $\g \rightarrow T_x\OO$ given by $X \mapsto X^*_x$. Hence,
the Jacobi identity always holds in this case.

$\bullet$ $u_1, u_2 \in T_x\OO$ and $u_3 \in T_x\OO^\perp$. In this case
we can write $u_1 = X^*_x$ and $u_2 = Y^*_x$ for some $X,Y \in \g$. Then, by
the definition of $[\cdot,\cdot]_0$:
\begin{align*}
	[[X^*_x,Y^*_x]_0,u_3]_0 &= [[X,Y]^*_x,u_3]_0  = [X,Y](u_3) \\
		&= X(Y(u_3)) - Y(X(u_3))  \\
		&= X([Y^*_x,u_3]_0) -Y([X^*_x,u_3]_0)  \\
		&= [X^*_x,[Y^*_x,u_3]_0]_0
		- [Y^*_x,[X^*_x,u_3]_0]_0
\end{align*}
which proves that the Jacobi identity holds in this case in general.

$\bullet$ $u_1 \in T_x\OO$ and $u_2, u_3 \in T_x\OO^\perp$. We can now
choose $X \in \g$ such that $u_1 = X^*_x$. Hence, using from
Lemma~\ref{lemma-omega-Omega} the fact that $\Omega_x$ is a homomorphism of
$\g$-modules, we obtain:
\begin{align*}
	[X^*_x,[u_2,u_3]_0]_0 &= [X^*_x,\Omega_x(u_2,u_3)^*_x]_0 
			= [X,\Omega_x(u_2,u_3)]^*_x  \\
		&= \Omega_x(X(u_2),u_3)^*_x + \Omega_x(u_2,X(u_3))^*_x  \\
		&= \Omega_x([X^*_x,u_2]_0,u_3)^*_x
		+ \Omega_x(u_2,[X^*_x,u_3]_0)^*_x  \\
		&= [[X^*_x,u_2]_0,u_3]_0 + [u_2,[X^*_x,u_3]_0]_0,
\end{align*}
which yields again the Jacobi identity without extra conditions.

$\bullet$ $u_1, u_2, u_3 \in T_x\OO^\perp$. The definition of
$[\cdot,\cdot]_0$ yields now:
$$
	[[u_1,u_2]_0,u_3]_0 = [\Omega_x(u_1,u_2)^*_x,u_3]_0 =
			\Omega_x(u_1,u_2)(u_3),
$$
and so the Jacobi identity is satisfied in this case exactly when $\Omega_x$
satisfies condition (3).

The above proves the equivalence in the statement, and so it remains to obtain
the last claim. For this we need to show that:
$$
	X([u_1,u_2]_0) = [X(u_1),u_2]_0	+ [u_1,X(u_2)]_0 
$$
for every $X \in \g$ and $u_1,u_2 \in T_x M$. This is now dealt with
through the following cases which just basically apply the definitions involved
and properties already considered.

$\bullet$ $u_1,u_2 \in T_x\OO$. Then, we can write $u_1 = Y^*_x, u_2 = Z^*_x$
for some $Y,Z \in \g$ and:
\begin{align*}
	X([Y^*_x,Z^*_x]_0) &= X([Y,Z]^*_x) = [\rho_x(X),[Y,Z]^*]_x
			= [X,[Y,Z]]^*_x  \\
		&= [[X,Y],Z]^*_x +[Y,[X,Z]]^*_x  \\
		&= [[X,Y]^*_x,Z^*_x]_0 +[Y^*_x,[X,Z]^*_x]_0  \\
		&= [[\rho_x(X),Y^*]_x,Z^*_x]_0 
			+[Y^*_x,[\rho_x(X),Z^*]_x]_0  \\
		&= [X(Y^*_x),Z^*_x]_0 + [Y^*_x,X(Z^*_x)]_0
\end{align*}

$\bullet$ $u_1 \in T_x\OO$ and $u_2 \in T_x\OO^\perp$. Now we can write
$u_1 = Y^*_x$ for $Y \in \g$ and:
\begin{align*}
	X([Y^*_x,u_2]_0) &= X(Y(u_2)) = [X,Y](u_2) + Y(X(u_2))  \\
		&= [[X,Y]^*_x,u_2]_0 + [Y^*_x,X(u_2)]_0  \\
		&= [[\rho_x(X),Y^*]_x,u_2]_0 + [Y^*_x,X(u_2)]_0  \\
		&= [X(Y^*_x),u_2]_0 + [Y^*_x,X(u_2)]_0  
\end{align*}

$\bullet$ $u_1,u_2 \in T_x\OO^\perp$. We now have:
\begin{align*}
	X([u_1,u_2]_0) &= X(\Omega_x(u_1,u_2)^*_x)
			= [\rho_x(X),\Omega_x(u_1,u_2)^*]_x  \\
		&= [X,\Omega_x(u_1,u_2)]^*_x  \\
		&= \Omega_x(X(u_1),u_2)^*_x + \Omega_x(u_1,X(u_2))^*_x  \\
		&= [X(u_1),u_2]_0 + [u_1,X(u_2)]_0
\end{align*}
\end{proof}

The following result will allow us to prove
Theorem~\ref{theorem-alghull-TOperp-irred}. It also provides an explicit
description of the Lie algebra structure considered in
Theorem~\ref{theorem-alghull-TOperp-irred}.

\begin{theorem}\label{theorem-alghull-TOperp}
Suppose that $G$ has finite center and real rank at least $2$, and that the
$G$-action on $M$ preserves a finite volume complete pseudo-Riemannian metric.
Also assume that $G$ acts ergodically on $M$ and that the foliation by
$G$-orbits is nondegenerate. 
Denote with $L$ the algebraic hull for the $G$-action on the bundle
$L(T\OO^\perp)$ and with $\mathfrak{l}$ its Lie algebra. In particular, there is
an embedding of Lie algebras $\mathfrak{l} \hookrightarrow \so(p,q)$, where
$(p,q)$ is the signature of the metric of $M$ restricted to $T\OO^\perp$.
If this embedding is surjective, then one of the following occurs:
\begin{enumerate}
\item the conclusion of Theorem~\ref{thm-int-normal} holds, or

\item $G$ is locally isomorphic to $\SO_0(p,q)$, $\dim(M) = (p+q)(p+q+1)/2$, 
for some $x \in M$ the bilinear map $\Omega_x$ is nonzero and
of curvature type and the Lie algebra structure on $T_x M$
obtained from Lemma~\ref{lemma-Omega-curvature} is isomorphic to either
$\so(p,q+1)$ or $\so(p+1,q)$.
\end{enumerate}
\end{theorem}
\begin{proof}
Since the $G$-orbits are nondegenerate, $G$ preserves a pseudo-Riemannian metric
on $T\OO^\perp$. Hence, the algebraic hull of $L(T\OO^\perp)$ for the
$G$-action can be embedded into the structure group $\mathrm{O}(p,q)$ for such
metric, where $(p,q)$ is the signature of $T\OO^\perp$. This yields the first
claim.

Let us now assume that the induced embedding $\mathfrak{l} \hookrightarrow
\so(p,q)$ is surjective. This implies that $L$ is a finite index subgroup of
$\mathrm{O}(p,q)$. By Section~3 of \cite{Zimmer-alghull}, we can write $L =
ZS$ where $S$ is semisimple without compact factors, $Z$ is compact and
centralizes $S$ and the product is almost direct. In particular, we
have a direct product $\mathfrak{l} = \mathfrak{z} \times \s$. Note that
for $p+q \leq 2$ the conclusion of Theorem~\ref{thm-int-normal} holds by
Corollary~\ref{cor-int-normal-rank}(1). Hence, we can assume from now on
that $p+q \geq 3$. Then, the only cases in which $\so(p,q)$ is not
simple is for $\so(4) \cong \so(3)\times\so(3)$ and $\so(2,2) \cong
\slinear(2,\R)\times \slinear(2,\R)$. And so, one of the following holds:
\begin{itemize}
\item $\s = 0$ and $\mathfrak{z} \cong \so(p,q)$ is compact, or
\item $\mathfrak{z} = 0$ and $\s \cong \so(p,q)$ is noncompact.
\end{itemize}
In the first case, we conclude that $\so(T_x\OO^\perp)$ is compact for every $x
\in M$, thus implying that the conclusion of
Theorem~\ref{thm-int-normal}
follows by Corollary~\ref{cor-int-normal-rank}(1). In particular, we can assume
that $\mathfrak{l} = \s \cong \so(p,q)$. Also, by the arguments in Section~3 of
\cite{Zimmer-alghull}, there is a surjection $\g \rightarrow \s$ of Lie algebras
which implies that $\g \cong \so(p,q)$.
Hence, $G$ is locally isomorphic to $\SO_0(p,q)$, and because of the
decomposition $TM = T\OO \oplus T\OO^\perp$ we also have $\dim(M) =
\dim(\SO_0(p,q)) + \dim(T_x\OO^\perp) = \dim(\SO_0(p,q)) + p + q =
(p+q)(p+q+1)/2$.

On the other hand, for $S \subset M$ given as in
Proposition~\ref{prop-g(x)}, if $\Omega|_S = 0$, then the conclusion of
Theorem~\ref{thm-int-normal} holds by
Lemma~\ref{lemma-omega-Omega}(2). Then, we can choose $x \in S$ such that
$\Omega_x \neq 0$. In particular, by Lemma~\ref{lemma-omega-Omega}(1) the
$\g$-module $T_x\OO^\perp$ is a nontrivial one. 
If we consider the representation
$\g \rightarrow \so(T_x\OO^\perp)$ that defines the $\g$-module structure of
$T_x\OO^\perp$ from Proposition~\ref{prop-g(x)}, then the above remarks show
that this representation is an isomorphism. We conclude that $T_x\OO^\perp$ is
a $\g$-module isomorphic to $\R^{p,q}$ with respect to the isomorphism $\g
\rightarrow \so(T_x\OO^\perp)$ thus obtained.

For $\left<\cdot,\cdot\right>_{p,q}$ the metric in $\R^{p,q}$ preserved by
$\so(p,q)$ we have the following.

\emph{Claim:} For every $c \in \R$, the map $T_c : \wedge^2\R^{p,q} \rightarrow
\so(p,q)$ given by:
$$
	T_c(u\wedge v) = c\left<v,\cdot\right>_{p,q}u -
			c\left<u,\cdot\right>_{p,q}v,
$$
where $u,v \in \R^{p,q}$, is a well defined homomorphism of
$\so(p,q)$-modules. Moreover, these maps exhaust all the $\so(p,q)$-module
homomorphisms $\wedge^2\R^{p,q} \rightarrow \so(p,q)$.

The only nontrivial part of the claim is the last statement, which we will now
prove. Let  $T : \wedge^2\R^{p,q} \rightarrow \so(p,q)$ be a homomorphism of
$\so(p,q)$-modules. Consider the map $T\circ T_1^{-1}$ which is a homomorphism
of $\so(p,q)$-modules $\so(p,q) \rightarrow \so(p,q)$. Hence, for $n=p+q$ and by
complexifying, the map $T\circ T_1^{-1}$ yields a homomorphism $\widehat{T} :
\so(n,\C) \rightarrow \so(n,\C)$ of $\so(n,\C)$-modules such that
$\widehat{T}|_{\so(p,q)} = T\circ T_1^{-1}$. In particular,
$\widehat{T}(\so(p,q))\subset \so(p,q)$ and so $\widehat{T}$ commutes with the
conjugation $\sigma$ of $\so(n,\C)$ whose fixed point set is $\so(p,q)$,
i.e.~$\widehat{T}\circ \sigma = \sigma\circ \widehat{T}$. Note that since
$\g\cong\so(p,q)$ is simple with real rank at least $2$, the Lie algebra
$\so(n,\C)$ is simple as well. Hence, the irreducibility of $\so(n,\C)$ as
$\so(n,\C)$-module implies, by Schur's Lemma, that there is a complex number $c$
such that $\widehat{T} = c\mathrm{Id}_{\so(n,\C)}$. But then, the condition
$\widehat{T}\circ \sigma = \sigma\circ \widehat{T}$ implies that $c$ is real and
so $T = T_c$.

By the Claim and Lemma~\ref{lemma-omega-Omega}(1), with respect to the above
established isomorphisms $\g \rightarrow \so(T_x\OO^\perp)$ and $T_x\OO^\perp
\cong \R^{p,q}$, we conclude that the map $\Omega_x$
corresponds to a homomorphism $T_c$ for some real $c \neq 0$. It is
straightforward to check that $T_c$ is of curvature type; in fact, it yields the
curvature of the pseudo-Riemannian manifolds of constant sectional curvature $c$
and signature $(p,q)$ (see Corollary~43 in page~80 of \cite{ONeill-book}).
Hence, Lemma~\ref{lemma-Omega-curvature} applies and it provides a Lie algebra
structure on $T_x M$. It remains to show that such structure is
isomorphic to either $\so(p,q+1)$ or $\so(p+1,q)$.

By our identifications, the Lie algebra structure on $T_x M$ is
isomorphic to the one obtained on $\so(p,q)\oplus\R^{p,q}$ from the map $T_c$
with the Lie brackets $[\cdot,\cdot]_c$ defined with formulas similar to those
in Lemma~\ref{lemma-Omega-curvature} replacing $\Omega_x$ with $T_c$. Also, it
is straightforward to show that the map given by:
\begin{align*}
	(\so(p,q)\oplus\R^{p,q},[\cdot,\cdot]_c) &\rightarrow
				\so(\R^{p+q+1},I_{p,q}(c)) \\
	X + u &\rightarrow \left(
			\begin{matrix}
			X & cu \\
			u^tI_{p,q} & 0
			\end{matrix}
			\right),
\end{align*}
is an isomorphism of Lie algebras, where:
$$
	I_{p,q}(c) = \left(
		\begin{matrix}
		I_{p,q} & 0 \\
		0 & c
		\end{matrix}
		\right).
$$
Since $\so(\R^{p+q+1},I_{p,q}(c))$ is clearly isomorphic to either $\so(p,q+1)$
or $\so(p+1,q)$ this yields (2) from our statement.
\end{proof}

We can now complete the following proof.

\begin{proof}[Proof of Theorem~\ref{theorem-alghull-TOperp-irred}]
Since the hypotheses of Theorem~\ref{theorem-alghull-TOperp} are a subset
of those of Theorem~\ref{theorem-alghull-TOperp-irred}, the conclusions of the
former hold. Since $M$ is weakly irreducible its universal covering space
cannot split isometrically and so part (2) of
Theorem~\ref{theorem-alghull-TOperp} necessarily holds at some $x \in
M$. By the definition of the Lie algebra structure in
$T_x M$ from Lemma~\ref{lemma-Omega-curvature} we conclude that (1) and (2) of
Theorem~\ref{theorem-alghull-TOperp-irred} are satisfied. 

Finally we observe that, by the proof of Theorem~\ref{theorem-alghull-TOperp},
we can assume that $x \in S$ where $S \subset M$ is given by
Proposition~\ref{prop-g(x)}. Consider the Lie algebra $\g(x)$ of local Killing
fields provided by Proposition~\ref{prop-g(x)}. Then, the last claim of
Lemma~\ref{lemma-Omega-curvature} and the definition of the representations
involved implies the last claim of Theorem~\ref{theorem-alghull-TOperp-irred}
for the Lie algebra $\g(x)$.
\end{proof}


\begin{thebibliography}{00}

\bibitem{GCT}
A.~Candel and R.~Quiroga--Barranco, Gromov's centralizer theorem,
Geom.~Dedicata {\bf 100} (2003), 123--155.

\bibitem{Gromov}
M.~Gromov, Rigid transformations groups, in G\'eom\'etrie diff\'erentielle,
Colloque G\'eom\'etrie et Physique de 1986  en l'honneur de Andr\'e Lichnerowicz
(D. Bernard and Y. Choquet-Bruhat, eds.), Hermann, 1988, 65--139.

\bibitem{Hernandez}
F.~G.~Hern\'andez-Zamora, Isometric splitting for actions of simple Lie
groups on pseudo-Riemannian manifolds,  Geom. Dedicata  {\bf 109}  (2004),
147--163.

\bibitem{Nadine-Annals}
N.~Kowalsky, Noncompact simple automorphism groups of Lorentz manifolds and
other geometric manifolds, Ann. of Math. (2) {\bf 144} (1996), no. 3, 611--640. 

\bibitem{Molino}
P.~Molino, Riemannian foliations, Translated from the French by Grant Cairns.
With appendices by Cairns, Y. Carri‚Äö√Ñ√∂‚àö‚Ä†‚àö‚àÇ¬¨¬®‚àö√úre, \'E. Ghys, E. Salem and V. Sergiescu.
Progress in Mathematics, 73. Birkh\"auser Boston, Inc., Boston, MA, 1988.

\bibitem{Nomizu}
K.~Nomizu, On local and global existence of Killing vector fields, Ann.~of
Math. {\bf 72} (1960), 105--120.

\bibitem{Gestur-Quiroga}
G.~Olafsson and R.~Quiroga-Barranco, Isometric actions of
$\SO_0(p,q)$ on low dimensional weakly irreducible pseudo-Riemannian
manifolds, In preparation.

\bibitem{ONeill-submersion}
B.~O'Neill, The fundamental equations of a submersion, Michigan Math. J. 
{\bf 13} (1966), 459--469.

\bibitem{ONeill-book}
B.~O'Neill, Semi-Riemannian geometry. With applications to relativity. Pure
and Applied Mathematics, 103. Academic Press, Inc., New York, 1983.

\bibitem{Quiroga-Annals}
R.~Quiroga-Barranco, Isometric actions of simple Lie groups on
pseudo-Riemannian manifolds, Ann. of Math. (2) {\bf 164} (2006), no. 6,
941--969.

\bibitem{Quiroga-AJM}
R.~Quiroga-Barranco, Arithmeticity of totally geodesic Lie foliations with
locally symmetric leaves, To appear in Asian J. Math.

\bibitem{Singer}
I.~M.~Singer, Infinitesimally homogeneous spaces, Comm. Pure Appl. Math. {\bf
13} (1960) 685--697. 

\bibitem{Szaro}
J.~Szaro, Isotropy of semisimple group actions on manifolds
with geometric structure, Amer. J. Math. {\bf 120} (1998), 129--158.

\bibitem{Zeghib}
A.~Zeghib, On affine actions of Lie groups, Math. Z. {\bf 227} (1998), no. 2,
245--262.

\bibitem{Zimmer-prog}
R.~J.~Zimmer, Actions of semisimple groups and discrete subgroups.
Proceedings of the International Congress of Mathematicians, Vol.
1, 2 (Berkeley, Calif., 1986), 1247--1258, Amer. Math. Soc.,
Providence, RI, 1987.

\bibitem{Zimmer-Lie}
R.~J.~Zimmer, Arithmeticity of holonomy groups of Lie foliations, J.
Amer. Math. Soc. {\bf 1} (1988), no. 1, 35--58.

\bibitem{Zimmer-alghull}
R.~J.~Zimmer, On the algebraic hull of an automorphism group of a principal
bundle, Comment. Math. Helv. {\bf 65} (1990), no. 3, 375--387. 

\bibitem{Zimmer-geometric}
R.~J.~Zimmer, Automorphism groups and fundamental groups of geometric
manifolds, in Differential Geometry: Riemannian Geometry (Los Angeles, CA,
1990), 693--710, Proc. Sympos. Pure Math. 54, Part 3, A. M. S., Providence, RI,
1993.


\end{thebibliography}
\end{document}